\documentclass[11pt]{article}
\usepackage{amsmath,amsthm,amssymb,amsfonts,fullpage,verbatim,bm,graphicx,enumerate,epstopdf,lscape,multirow,adjustbox,changepage,subfig}
\usepackage[dvipsnames,usenames]{color}
\usepackage[pdftex,colorlinks=true,pdfpagemode=none,linkcolor=blue,
   citecolor=blue,plainpages=true]{hyperref}
\usepackage{tikz,verbatim,bibentry}
\usepackage{bm,mathrsfs}
\usepackage[comma,authoryear]{natbib}

\newlength{\offsetpage}
\newenvironment{widepage}{\begin{adjustwidth}{-\offsetpage}{-\offsetpage}%
    \addtolength{\textwidth}{2\offsetpage}}%
{\end{adjustwidth}}

\newtheorem{theorem}{Theorem}
\newtheorem{proposition}{Proposition}

\theoremstyle{definition}
\newtheorem{rmk}{Remark}
\newtheorem{exmp}{Example}

\newtheorem{assumption}{Assumption}

\def\Cov{\textsf{Cov}}

\def\E{\textsf{E}}
\def\P{\textsf{P}}
\def\tr{\mathrm{tr}}

\begin{document}
\title{High-dimensional Linear Regression for Dependent Data with Applications to Nowcasting}
\author{Yuefeng Han and Ruey S. Tsay}
\date{University of Chicago}
\maketitle

\begin{abstract}
\normalsize
\indent Recent research has focused on $\ell_1$ penalized least squares (Lasso) estimators for high-dimensional linear regressions in which the number of covariates $p$ is considerably larger than the sample size $n$. However, few studies have examined the properties of the estimators when the errors and/or the covariates are serially dependent.
In this study, we investigate the theoretical properties of the Lasso estimator for a linear regression with a random design and weak sparsity under serially dependent and/or nonsubGaussian errors and covariates. In contrast to the traditional case, in which the errors are independent and identically distributed and have finite exponential moments, we show that $p$ can be at most a power of $n$ if the errors have only finite polynomial moments. In addition,
the rate of convergence becomes slower owing to the serial dependence in the errors and the covariates. We also consider the sign consistency of the model selection using the Lasso estimator when there are serial correlations in the errors or the covariates, or both. Adopting the framework of a functional dependence measure, we describe how the rates of convergence and the selection consistency of the estimators depend on the dependence measures and moment conditions of the errors and the covariates. Simulation results show that a Lasso regression can be significantly more powerful than a mixed-frequency data sampling regression (MIDAS) and a Dantzig selector in the presence of irrelevant variables. We apply the results obtained for the Lasso method to nowcasting with mixed-frequency data, in which serially correlated errors and a large number of covariates are common. The empirical results show that the Lasso procedure outperforms the MIDAS regression and the autoregressive model with exogenous variables in terms of both forecasting and nowcasting.

\vspace{0.3in}
MSC2000 Subject Classification: Primary 62J05; Secondary 62M10.

\end{abstract}
\newpage

\section{Introduction}
The past two decades have witnessed significant developments in high-dimensional linear regression analyses. Consider the following linear regression for the
response variable $y_i$ and the covariate vector ${\bf x}_i$:
\begin{eqnarray} \label{model}
y_i={\bf x}_i^T\beta+e_i,\quad 1\leq i\leq n,
\end{eqnarray}
where $\beta\in \mathbb R^p$ consists of unknown coefficients, $e_i$ is an error term, and ${\bf x}_i^T$ denotes the transpose of  the covariate vector ${\bf x}_i$.
Denote the dimension of ${\bf x}_i$ by $p$. In matrix form, we can write the model as $Y=X\beta+e$, where $Y$ is an $n\times 1$ response vector, $X$ is an $n\times p$ design matrix, and $e$ is an $n\times 1$ error vector. Under certain sparsity conditions on $\beta$, many studies have focused on the $\ell_1$ penalized least squares (Lasso) estimator of $\beta$ when the number of variables $p$ can be much larger than the sample size $n$; see \citet{efron2004}, \citet{zhao2006}, and \citet{meinshausen2009}, among others. Other related approaches include the Dantzig selector of \citet{candes2007}, adaptive Lasso of \citet{zou2006}, group Lasso
by \citet{yuan2006}, and SCAD estimator of \citet{fan2001}, among others. The theoretical properties of those estimators have been established in the literature under the independence assumption; see, for example, \citet{bickel2009} and \citet{buhlmann2011}. Here, we focus on the Lasso estimator  defined as
\begin{eqnarray} \label{las}
\hat\beta=\arg\min_{\beta}\left(\frac{1}{2} |Y-X\beta|_2^2+\lambda |\beta|_1\right),
\end{eqnarray}
where $\lambda\geq0$ is a tuning parameter that controls the level of sparsity
in $\hat\beta$.

\indent Much of the available research dedicated to the Lasso problem examines the case of large $p$ and small $n$ when the design matrix is static and the errors are independent and identically distributed (i.i.d.) random variables. On the other hand, in many real applications, ${\bf x}_i$ consists of stochastic random variables that might be dynamically dependent, or $e_i$ is serially dependent, or both. Despite considerable recent work on Lasso estimators, few studies examine the theoretical properties of the estimates when the observations are dependent. \citet{wang2007} proposed a Lasso estimator for a regression model with autoregressive errors. \citet{Shuva2012} investigated the Lasso estimator for weakly dependent errors. Both studies concentrate on the case when $n$ is greater than $p$. More recently, \citet{basu2015} investigated the theoretical properties of Lasso estimators using a random design for high-dimensional Gaussian processes. \citet{kock2015} established the oracle inequalities of the Lasso for Gaussian errors in stationary vector autoregressive models. \citet{wu2016} analyzed the Lasso estimator with a fixed design matrix, and assumed that a restricted eigenvalue condition is satisfied. \citet{Medeiros2016} studied the asymptotic properties of the adaptive Lasso when the errors are nonGaussian and may be conditionally heteroskedastic. The goal of this study is to investigate the limiting properties of Lasso estimators for Model (\ref{model}) in the presence of serial dependence in both the covariate vector ${\bf x}_i$ and the errors.
We establish the rate of convergence and provide the sign consistency of the Lasso estimator under the weak sparsity condition. Our results extend beyond those of a fixed design and exact sparsity time series; thus, we do not assume a restricted eigenvalue condition on
either the sample or on the population covariance matrix.

In practice, many important macroeconomic variables are not sampled at the same frequency. For example, gross domestic product (GDP) data are available quarterly, industrial production data are published monthly, and most interest rate data are available daily. Analyzing such data jointly is referred to as a mixed-frequency data analysis. In the econometrics literature, \citet{ghysels2004} proposed a mixed-data sampling (MIDAS) approach to analyze such data. In particular, they use newly available high-frequency data to improve the prediction of a lower-frequency macroeconomic variable of interest, and refer to such predictions as {\em nowcasting}. Consider, for example, the problem of predicting the quarterly GDP growth rate $y_{n+1}$ at the forecast origin $i=n$. Here, the time interval is a quarter. Traditional forecasting methods employ quarterly
data available at $i=n$ to build a model, after which, they use the fitted model to perform a prediction. In practice, some monthly and daily data become available during the quarter $i=n+1$. Nowcasting uses newly available monthly and daily data to update its prediction of $y_{n+1}$. Therefore, the term nowcasting means taking advantage of high-frequency data within a given quarter to update the predictions of GDP growth rate of that quarter. In short, the basic principle of nowcasting is the exploitation of information published at higher frequencies than the target variable of interest in order to obtain an improved prediction  before the official lower-frequency data become available. Because high-frequency data are relatively common in practice, employing many covariates is common in nowcasting. Therefore,
Model (\ref{model}), with dependent covariates and errors, is applicable to nowcasting, and the Lasso method is highly relevant.
The mixed-data sampling approach of \citet{ghysels2004} has proven useful for various forecasting and nowcasting purposes. We compare the performance of the Lasso regression with that of the MIDAS regression and the
autoregressive model with exogenous variables (ARX).
To the best of our knowledge, this is the first study
to apply a Lasso regression to nowcasting. Simulation studies and empirical studies show that the
Lasso estimator outperforms the existing MIDAS regression and ARX model.

The rest of the paper is organized as follows. Section 2 defines the
high-dimensional dependence measure, adopting the concept of \citet{wu2005}.
Section 3 deals with rates of convergence of Lasso estimators.
The model selection consistency of Lasso estimators is given in Section 4, and simulation studies are carried out in Section 5. Section 6 considers real-data examples,
 including forecasting and nowcasting applications.

\indent We begin with 
some basic definitions. Throughout the paper, for a matrix $A=(a_{ij})\in\mathbb R^{p\times p}$, define the spectral norm $\rho(A)=\sup_{|x|\leq 1}|Ax|_2$,
the Frobenius norm $|A|_F = (\sum_{ij}a_{ij}^2)^{1/2}$, and the infinity norm $|A|_\infty = \max_{1\leq i,j\leq p} |a_{ij}|$. For a vector ${a}=(a_1,...,a_p)^T\in\mathbb R^p$, define the vector $q$ norm ${|a|_q=(\sum_{i=1}^p |a_i|^q)^{1/q}}$, for $1\leq q < \infty$. Let $|a|_\infty=\max_{1\leq i\leq p}|a_i|$ and $|a|_0=\#\{i: a_i\neq 0\}$.
For a random variable $\xi\in \mathcal{L}^k $, denote the $q$-norm by $\|\xi\|_q = (\E |\xi|^q)^{1/q}$, for $1\leq q \leq k$. For two sequences of real numbers $\{a_n\}$ and $\{b_n\}$, write $a_n=O(b_n)$ if there exists a constant $C$ such that $|a_n|\leq C |b_n|$ holds for all sufficiently large $n$, write $a_n=o(b_n)$ if $\lim_{n\to\infty} a_n/b_n =0$, and write $a_n\asymp b_n$ if there are positive constants c and C, such that $c\leq a_n/b_n\leq C$ for all sufficiently large $n$. Denote $a\wedge b=\min\{a,b\}$ and $a\vee b=\max\{a,b\}$.

\section{High-Dimensional Time Series}
Let $\varepsilon_i$, for $i\in\mathbb Z$, be i.i.d. random vectors and the $\sigma$-field $\mathcal F_i=(\cdots,\varepsilon_{i-1},\varepsilon_i)$. In our random-design setting, we assume that in Model (\ref{model}), the covariate process (${\bf x}_i,i=1,...,n$) is high-dimensional and weakly stationary, and of the form
\begin{eqnarray}\label{eq3}
{\bf x}_i=(g_1(\mathcal F_i),...,g_p(\mathcal F_i))^T,
\end{eqnarray}
and the error $e_i$ satisfies
\begin{eqnarray}\label{eq4}
e_i=g_e(\mathcal F_i),
\end{eqnarray}
where $g_1(\cdot), \ldots, g_p(\cdot)$ and $g_e(\cdot)$ are measurable functions in $\mathbb R$, such that ${\bf x}_i$ is well defined. In the scalar case with $p=1$, (\ref{eq3}) and (\ref{eq4}) include a very general class of stationary processes (see \citet{wiener1958}, \citet{rosenblatt1971}, \citet{priestley1988}, \citet{tong1990}, \citet{tsay2005}, \citet{wu2005}). They also allow models with homogeneous or heteroscedastic errors; see Example 1 of Section 3. In the homogeneous case, the covariate process (${\bf x}_i$) and the errors ($e_i$) can be independent of each other. \\
\indent Following \citet{wu2005}, we define the functional dependence measure
\begin{eqnarray}
\delta_{i,q,j}&=&\|x_{ij}-x_{ij}^*\|_q 
=\|g_j(\mathcal F_i)-g_j(\mathcal F_i^*)\|_q, \\
\delta_{i,q,e}&=&\|e_i-e_i^*\|_q 
=\|g_e(\mathcal F_i)-g_e(\mathcal F_i^*)\|_q,
\end{eqnarray}
where the coupled process $x_{ij}^*=g_j(\mathcal F_i^*)$ and $e_i^*=g_e(\mathcal F_i^*)$. Here $\mathcal F_i^*=(...,\varepsilon_{-1},\varepsilon_0',\varepsilon_1,...,\varepsilon_{i-1},\varepsilon_i)$ and $\varepsilon_0',\varepsilon_l$, for $l\in\mathbb Z$, are i.i.d. random variables. We assume short-range dependence, such that
\begin{eqnarray}
\Delta_{m,q,j}:=\sum_{i=m}^{\infty} \delta_{i,q,j} <\infty, \\
\Delta_{m,q,e}:=\sum_{i=m}^{\infty} \delta_{i,q,e} <\infty.
\end{eqnarray}
Then, for fixed $m$, $\Delta_{m,q,j}$, and $\Delta_{m,q,e}$, measure the cumulative effect of $\varepsilon_0$ on $(x_{ij})_{i\geq m}$ and $(e_i)_{i\geq m}$. We introduce the following dependence-adjusted norm (DAN):
\begin{eqnarray}
\|x_{.j}\|_{q,\alpha}&=&\sup_{m\geq0}(m+1)^\alpha\Delta_{m,q,j},\quad \alpha\geq 0. \\
\|e_.\|_{q,\alpha}&=&\sup_{m\geq0}(m+1)^\alpha\Delta_{m,q,e},\quad \alpha\geq 0.
\end{eqnarray}
It can happen that, owing to the dependence, $\|e_.\|_{q,\alpha}=\infty$, while $\|e_i\|_{q}<\infty$. Because $e_0=\sum_{l=-\infty}^0(\E(e_0|\mathcal F_{l})-\E(e_0|\mathcal F_{l-1}))$, we have
\begin{eqnarray}
\|e_0\|_q \leq \sum_{l=0}^\infty \|\E(e_0|\mathcal F_{-l})-\E(e_0|\mathcal F_{-l-1})\|_q=\sum_{l=0}^\infty \|\E(e_l-e_l^*|\mathcal F_0)\|_q \leq\sum_{l=0}^\infty\|e_l-e_l^*\|_q = \|e_.\|_{q,0},
\end{eqnarray}
by stationarity. If $e_i$, for $i\in\mathbb Z$, are i.i.d., the DAN $\|e_.\|_{q,\alpha}$ and the $\mathcal L^q$ norm $\|e_0\|_q$ are equivalent, in the sense that $\|e_0\|_q\leq \|e_.\|_{q,\alpha}\leq 2\|e_0\|_q$. \\
\indent To account for the cross-sectional dependence of the $p$-dimensional stationary process (${\bf x}_i$), we define the $\mathcal L^\infty$ functional dependence measure and its corresponding DAN (see \citet{chen2013}, \citet{danna2017}), as follows:
\begin{eqnarray*}
&&\omega_{i,q}=\| \max_{1\leq j\leq p} |x_{ij}-x_{ij}^*| \|_q,\\
&&\| |{\bf x}_{.}|_\infty\|_{q,\alpha}=\sup_{m\geq0}(m+1)^\alpha\Omega_{m,q},\quad \alpha\geq 0,\quad\text{and } \Omega_{m,q}=\sum_{i=m}^\infty \omega_{i,q}.
\end{eqnarray*}
Additionally, we define
\begin{eqnarray*}
\Psi_{q,\alpha}=\max_{1\leq j\leq p}\|x_{.j}\|_{q,\alpha}\quad\text{and }\quad \Upsilon_{q,\alpha} = \left(\sum_{j=1}^p \|x_{.j}\|_{q,\alpha}^q \right)^{1/q},
\end{eqnarray*}
where $\Psi_{q,\alpha}$ and $\Upsilon_{q,\alpha}$ can be viewed as the uniform and the overall DANs of (${\bf x}_i$), respectively. Clearly, $\Psi_{q,\alpha} \leq \| |{\bf x}_{.}|_\infty\|_{q,\alpha}\leq \Upsilon_{q,\alpha}$.

Next, we provide an example of high-dimensional time series to illustrate the univariate and multivariate DAN scale.
\begin{exmp} \label{exp:vma}
Let $\varepsilon_{ij}$, for $i,j\in\mathbb Z$, be i.i.d. random variables with mean zero and variance one, and with finite $q$th moments, $q >2$. Furthermore, let $A_i$, for $i\geq0$, be $p\times d$ coefficient matrices with real entries, such that $\sum_{i=0}^\infty \tr(A_iA_i^T)<\infty$. Write $\varepsilon_i= (\varepsilon_{i1},...,\varepsilon_{id})^T$. Then, by Kolmogorov's three-series theorem, the linear process
\begin{eqnarray}\label{eq:vma}
{\bf x}_i=\sum_{l=0}^\infty A_l\varepsilon_{i-l}
\end{eqnarray}
exists. Denote $A_l=(a_{l;jk})_{1\leq j\leq p,1\leq k\leq d}$, and $A_{l,j.}$ is the $j$th row of $A_l$. By Burkholder's inequality, $\|A_{l,j.}\varepsilon_0\|_{q} \le \sqrt{q-1} |A_{l,j.}|_2 \|\varepsilon_{00}\|_{q}$. We assume that the linear process satisfies the decay condition
\begin{eqnarray}\label{eq:linear}
\max_{j\le p} |A_{l,j.}|_2\le K_1(1\vee l)^{-\theta},
\end{eqnarray}
for all $l\ge0$, where $\theta>1/2$ and $K_1>0$. If $\theta>1$, (\ref{eq:linear}) implies short-range dependence (SRD) because the auto-covariance matrices $\Sigma_k=\sum_{l=0}^\infty A_l A_{l+k}^T$ are absolutely summable. On the other hand, if $1>\theta>1/2$, then (${\bf x}_i$) in (\ref{eq:vma}) may not have summable auto-covariance matrices, thus allowing for long-range dependence (LRD). The classical literature on LRD focuses primarily on the univariate case, $p=1$. Then, under the SRD case, the DANs have the following bounds:
\begin{eqnarray}
&&\Psi_{q,\alpha}=\max_{1\leq j\leq p}\|x_{.j}\|_{q,\alpha} =  \max_j \sup_{m\ge0} (m+1)^{\alpha} \sum_{i=m}^\infty\| A_{i,j.}\varepsilon_0 \|_q\le K_1K_2 \|\varepsilon_{00}\|_{q}, \\
&&\| |{\bf x}_{.}|_\infty\|_{q,\alpha} = \sup_{m\ge0} (m+1)^{\alpha} \sum_{i=m}^\infty\| \max_j |A_{i,j.}\varepsilon_0| \|_q  \le K_1K_2 p^{1/q}\|\varepsilon_{00}\|_{q},
\end{eqnarray}
where $\alpha=\theta-1$ and the constant $K_2$ depends only on $\theta$ and $q$.
\end{exmp}

In this paper, we use the DANs $\| |{\bf x}_{.}|_\infty\|_{q,\alpha}$, $\Psi_{q,\alpha}$, and $\Upsilon_{q,\alpha}$ to study the limiting properties of Lasso estimators in the presence of serial dependence. These adjusted norms are more convenient
than the commonly used mixing conditions for handling serial dependence in high-dimensional
time series.

\section{Convergence Rate of the Lasso Estimator}
In this section, we present the main results on convergence rate of the Lasso estimator for dependent data. In the low-dimensional case, the consistency of $\hat\beta$ relies on the assumption that the sample covariance matrix converges to the population covariance matrix. In the high-dimensional case ($n\ll p$), it requires that $|X(\hat\beta-\beta)|_2$ is small only when $|\hat\beta-\beta|_2$ is small. Let
$\hat\Sigma=(\hat\sigma_{jk})_{1\leq j,k\leq p}
=n^{-1}\sum_{i=1}^n x_i x_i^T$ be the sample covariance. Typically, researchers assume with high probability that the following restricted strong convexity condition holds:
\begin{eqnarray}\label{cond:rsc}
u'\hat\Sigma u \ge \kappa_1 |u|_2^2- \kappa_2 g(n,p) |u|_1^2,
\end{eqnarray}
for all $u\in\mathbb R^p$, where $\kappa_1,\kappa_2$ are positive constants, and $g(n,p)$ is a function of the sample size $n$ and the ambient dimension $p$. This can be viewed as an analogous sufficient condition in the high-dimensional case. As shown in the proof of Theorem \ref{thm1}, the restricted strong convexity condition for the sample covariance matrix holds with high probability under certain conditions.


To establish our theoretical results, we first impose a weak sparsity condition.
\begin{assumption}\label{asmp:wsc}
There exists some $0\le \theta< 1$, with a uniform radius $K_\theta$, such that
\begin{eqnarray}\label{eq:wsc}
\sum_{j=1}^p |\beta_j|^\theta\le K_\theta.
\end{eqnarray}
\end{assumption}

The following theorem shows that the $L_2$ and $L_1$ convergence rates of $\hat\beta$ to $\beta$ depend on the moment condition and on the temporal and cross-sectional dependence conditions.

\begin{theorem}\label{thm1}
Denote the population covariance matrix by $\Sigma=(\sigma_{jk})
=[Cov(x_{ij},x_{ik})].$ Suppose the minimum eigenvalue of $\Sigma$ 
satisfies
$\lambda_{\min}(\Sigma)\ge \kappa>0.$
Assume that $\Psi_{\gamma,\alpha_X}=\max_j\|x_{.j}\|_{\gamma,\alpha_X}=M_X<\infty$, and $\|e_.\|_{q,\alpha_e}=M_e<\infty$, where $q>2,\gamma>4$ and $\alpha_X, \alpha_e>0$.
Define
\[ \nu = \left\{\begin{array}{ll} 1 & \mbox{if $\alpha_X \geq 1/2-2/\gamma$,}  \\ \gamma/4-\alpha_X\gamma/2 & \mbox{if $\alpha_X <1/2-2/\gamma$.}\end{array}\right. \]
 Assume $\tau=q\gamma/(q+\gamma)>2$ and let $\alpha=\min(\alpha_X,\alpha_e)$. Define
 \[ \rho= \left\{\begin{array}{ll} 1 &  \mbox{if $\alpha \geq 1/2-1/\tau$,} \\
  \tau/2-\alpha\tau &  \mbox{if $\alpha <1/2-1/\tau$.}\end{array}\right. \]
Denote $\omega=\sqrt{\log p/n} M_X^2 + n^{2\nu/\gamma-1} (\log p)^{3/2} \| |{\bf x}_{.}|_\infty\|_{\gamma,\alpha_X}^2$.
Suppose Assumption \ref{asmp:wsc} holds. Then, for any $\lambda$ such that
$$ \lambda\gtrsim
\sqrt{\log p/n}M_e M_X + n^{\rho/\tau-1}(\log p)^{3/2} M_e \| |{\bf x}_{.}|_\infty\|_{\gamma,\alpha_X}, $$
and $K_\theta \omega \lambda^{-\theta} \le C$ for some positive constant $C$, any Lasso solution $\hat\beta$ satisfies
\begin{eqnarray}
|\hat\beta-\beta|_2 &\lesssim& \sqrt{K_\theta}\left(\frac{\lambda}{\kappa}\right)^{1-\theta/2} \label{eqw1}, \\
|\hat\beta-\beta|_1 &\lesssim& K_\theta\left(\frac{\lambda}{\kappa}\right)^{1-\theta} \label{eqw2},
\end{eqnarray}
with probability at least $1-C_1 (\log p)^{-\gamma/2}-C_2 p^{-C_3}-C_4 (\log p)^{-\tau}$, where $C_1,..., C_4$ are positive constants.
\end{theorem}

In the special case $\theta=0$, the quantity of weak sparsity corresponds to an exact sparsity constraint; that is, $\beta$ has at most $s:=K_0$ nonzero entries. The following theorem shows the convergence rate of $\hat\beta$ 
and the prediction error $|X(\hat\beta-\beta)|_2^2$ for the exact sparsity case.

\begin{theorem}\label{thm1b}
Suppose the same conditions of Theorem \ref{thm1} hold. If $|\beta|_0=s$, $\kappa\asymp 1$, and
$$n \gtrsim M_X^4s^2\log p + s^{1/(1-2\nu/\gamma)} (\log p)^{3/(2-4\nu/\gamma)} \| |{\bf x}_{.}|_\infty\|_{\gamma,\alpha_X}^{2/(1-2\nu/\gamma)},$$
then, for any $\lambda$ such that
$$ \lambda\gtrsim
\sqrt{\log p/n}M_e M_X + n^{\rho/\tau-1}(\log p)^{3/2} M_e \| |{\bf x}_{.}|_\infty\|_{\gamma,\alpha_X}, $$
any Lasso solution $\hat\beta$ satisfies
\begin{eqnarray}
|\hat\beta-\beta|_2 &\lesssim& \lambda\sqrt{s}/\kappa \label{eqa}, \\
|\hat\beta-\beta|_1 &\lesssim& \lambda s/\kappa \label{eqb}, \\
|X(\hat\beta-\beta)|_2^2/n &\lesssim& \lambda^2 s/\kappa \label{eqc},
\end{eqnarray}
with probability at least
$1-C_1 (\log p)^{-\gamma/2}-C_2 p^{-C_3}-C_4 (\log p)^{-\tau}$.
\end{theorem}

\begin{rmk}
In the exact sparsity case, instead of the condition
$\lambda_{\min}(\Sigma)\ge \kappa>0,$ we may require that the restricted eigenvalue assumption RE($s$,3) of \citet{bickel2009} holds for the population covariance matrix
$\Sigma$; that is
\begin{eqnarray}\label{RE}
\kappa:=\min_{J\subseteq\{1,...,p\},|J|_0\leq s}\ \min_{u\neq0, |u_{J^c}|_1 \leq 3|u_{J}|_1} u'\Sigma u/|u|_2^2 >0,
\end{eqnarray}
where $J^c$ is the complement of the set $J$, that is, $J^c=\{1,2,...,p\} \backslash J$, and $u_J$ is defined as a modification of $u$ by setting its elements outside $J$ to zero.
All bounds (\ref{eqa}), (\ref{eqb}), and (\ref{eqc}) still hold with high probability.
\end{rmk}

\begin{rmk}
The best known convergence rate of Lasso estimators for i.i.d. subGaussian data requires that $K_\theta (\log p/n)^{1-\theta/2} \le C$, for some positive constant $C$.
Our theorems require that $K_\theta \omega \lambda^{-\theta} \le C$, where
$$\omega=\sqrt{\log p/n} M_X^2 + n^{2\nu/\gamma-1} (\log p)^{3/2} \| |{\bf x}_{.}|_\infty\|_{\gamma,\alpha_X}^2,$$
and
$$ \lambda\gtrsim
\sqrt{\log p/n}M_e M_X + n^{\rho/\tau-1}(\log p)^{3/2} M_e \| |{\bf x}_{.}|_\infty\|_{\gamma,\alpha_X}. $$
The second terms in $\omega$ and $\lambda$ are introduced by the heavy tails, and thus are unavoidable. In other words, under heavy-tailed distributions in some cases, the allowed dimension $p$ for Lasso methods can be at most a power of the sample size $n$.

In the exact sparsity case, we require $n \gtrsim M_X^4s^2\log p + s^{1/(1-2\nu/\gamma)} (\log p)^{3/(2-4\nu/\gamma)} \| |{\bf x}_{.}|_\infty\|_{\gamma,\alpha_X}^{2/(1-2\nu/\gamma)}$. One may argue that the first term $M_X^4s^2\log p$ can be further improved to $M_X^4s \log p$ for short-range temporal dependence data, in accordance with i.i.d. subGaussian data. However, we cannot achieve this because even the optimal Bernstein-type inequality for nonlinear weakly dependent data is still an open problem. The best known result is proposed by \citet{merlevede2009}.
\end{rmk}

\begin{rmk}
Based on Theorem \ref{thm1b}, we have the following cases.
Assume $M_X\asymp 1$ and $M_e\asymp 1$. Under the weak cross-sectional dependence $\| |\mathbf x_.|_\infty\|_{\gamma,\alpha_X} \asymp p^{1/\gamma}$, which holds if the $p$ components $x_{ij}$ ($1\leq j\leq p$) are nearly independent, the required sample size for exact sparsity is $n\gtrsim s^2\log p+ s^{1/(1-2\nu/\gamma)} (\log p)^{3/(2-4\nu/\gamma)} p^{2/(\gamma-2\nu)}$, and the regularization parameter satisfies $\lambda\gtrsim \sqrt{\log p/n} + n^{\rho/\tau-1}(\log p)^{3/2}p^{1/\gamma}$. By comparison, the Bonferroni inequality and Lemma 1 in the Appendix yield $n\gtrsim s^2\log p+ s^{1/(1-2\nu/\gamma)} p^{4/(\gamma-2\nu)}$ and $\lambda\gtrsim \sqrt{\log p/n} + n^{\rho/\tau-1}p^{1/\tau}$, respectively.

\indent In addition, under the strong cross-sectional
dependence $\| |\mathbf x_.|_\infty\|_{\gamma,\alpha_X} \asymp 1$, which holds if the $p$ components $x_{ij}$ ($1\leq j\leq p$) are linear combinations of fixed random variables, the required sample size for exact sparsity is $n\gtrsim s^2\log p+ s^{1/(1-2\nu/\gamma)} (\log p)^{3/(2-4\nu/\gamma)}$, and the regularization parameter satisfies  $\lambda\gtrsim \sqrt{\log p/n} + n^{\rho/\tau-1}(\log p)^{3/2}$.
\end{rmk}

\indent  Next, we apply the results of Theorem \ref{thm1} in an example.
\begin{exmp} \label{example1}
Consider the autoregressive model with exogenous variables, that is,
the ARX($a,b$) model:
\begin{eqnarray}
y_i=\sum_{l=1}^a\phi_{l}y_{i-l}+\sum_{l=0}^b\psi_{l}'{\bf z}_{i-l}+e_i=\beta'{\bf x}_i+e_i,
\end{eqnarray}
where $a$ and $b$ are nonnegative integers, $e_i$ follows a GARCH(1,1) model defined below,
and ${\bf z}_i$ is a linear process defined by
\begin{eqnarray}
{\bf z}_i=\sum_{l=0}^\infty A_l\varepsilon_{i-l},
\end{eqnarray}
where the random variables $\varepsilon_{ij}$ and coefficient matrices
$A_l$ are given in Example 1, with $E|\varepsilon_{ij}|^\gamma < \infty$
and $\gamma > 2$.
Assume the roots of the polynomial $1-\sum_{l=1}^a\phi_l B^l$ are
outside the unit circle, which ensures the stationarity of the autoregressive part of the model. In addition, assume the population covariance matrix $\Sigma=\E{\bf x}_i{\bf x}_i'$ is positive definite. \\
\indent Let
\begin{eqnarray}
e_i=\sqrt{h_i}\eta_i,\qquad h_i=\pi_0+\pi_1 e_{i-1}^2+\pi_2 h_{i-1},
\end{eqnarray}
with $\pi_0>0$, $\pi_1\ge 0$, $\pi_2\ge 0$, and $\E(\pi_1+\pi_2\eta_{i-1}^2)^{q/2}<\infty$, $q>4$. Then, it is easy to show that $\|e_.\|_{q,\alpha_e}<\infty$.

Again, by Burkholder's inequality, $\|A_{l,j.}\varepsilon_0\|_{\gamma} \le \sqrt{\gamma-1} |A_{l,j.}|_2 \|\varepsilon_{00}\|_{\gamma}$. If there exist constants $K_1>1$ and $\alpha_Z>0$, such that $\max_{j\le p} |A_{l,j.}|_2\le K_1(l+1)^{-1-\alpha_Z}$ holds for all $l\ge 0$, then we have $\max_j \|z_{.j}\|_{\gamma,\alpha_Z} \le K_1K_2 \|\varepsilon_{00}\|_{\gamma}$, where the constant $K_2$ depends only on $\alpha_Z$ and $\gamma$. Together with the assumption that the roots of the polynomial $1-\sum_{l=1}^a\phi_l  B^l$ are outside the unit circle, we ensure $\max_j \|x_{.j}\|_{\gamma,\alpha_Z}<\infty$.
\end{exmp}

\section{Model Selection Consistency}
In this section, we extend the asymptotic properties of the sign consistency for model selection, using the Lasso, to the dependent setting. The sign consistency of the Lasso was first introduced by \citet{zhao2006}. Without loss of generality, write $\beta=(\beta_1,...,\beta_s,...,\beta_p)'$, where $\beta_j\neq0$ if $j\le s$, and $\beta_j=0$ if $j>s$. That is, the first $s$ predictors are relevant variables. Denote $\beta=(\beta_{(1)}',\beta_{(2)}')'$, where $\beta_{(1)}$ is an $s\times 1$ vector. Correspondingly, for any $i$, denote ${\bf x}_i=({\bf x}_{i(1)}',{\bf x}_{i(2)}')'$ and $X=({\bf x}_1,...,{\bf x}_n)'=(X_{(1)},X_{(2)})$, where $X_{(1)}$ is an $n\times s$ sub-matrix of relevant variables, and $X_{(2)}$ is an $n\times (p-s)$ sub-matrix of irrelevant variables. Similarly, consider the partition of the covariance matrix as 
$$\Sigma=\left(\begin{matrix}
\Sigma_{11}&\Sigma_{12}\\
\Sigma_{21}&\Sigma_{22}
\end{matrix}\right),$$
where $\Sigma_{11}=\E {\bf x}_{i(1)} {\bf x}_{i(1)}'$ is an $s\times s$ sub-matrix
associated with the relevant variables.

\indent We impose the following assumptions.
\begin{assumption} \label{assump1}
For any $1\le i\le n$, $\E (x_{ik} | X_{(1)}, e)=\Sigma_{2k,1} \Sigma_{11}^{-1} x_{i(1)}$, where $\Sigma_{2k,1}$ is the $k$th row of $\Sigma_{21}$.
\end{assumption}
Define $z_{ik}=x_{ik}-\E( x_{ik}|X_{(1)}, e)$, for $s+1\le k\le p$, and ${\bf z}_i=(z_{i,s+1},...,z_{i,p})'$.
\begin{assumption} \label{assump2}
There exists $L>0$, such that $\min_{1\leq j\leq s}|\beta_j| \ge L$.
\end{assumption}
\begin{assumption} \label{assump3}
There exists a constant $N_1>0$, such that
$$\inf_{|\zeta|_2=1}\zeta'\Sigma_{11}\zeta =N_1.$$
\end{assumption}
\begin{assumption} \label{assump4}
There exists a positive constant $\eta\in(0,1)$, such that
\begin{eqnarray} \label{irrep}
|\Sigma_{21}\Sigma_{11}^{-1}\text{sign}(\beta_{(1)})|_\infty \leq 1-\eta.
\end{eqnarray}
\end{assumption}

Assumption \ref{assump1} explicitly defines how the irrelevant variables depend on the relevant variables and the errors. Note that $\Cov(\Sigma_{2k,1} \Sigma_{11}^{-1} x_{i(1)}, x_{ik}-\Sigma_{2k,1} \Sigma_{11}^{-1} x_{i(1)})=0$ always holds, for all $s+1\le k\le p$. That is, $\Sigma_{2k,1} \Sigma_{11}^{-1} x_{i(1)}$ and $ x_{ik}-\Sigma_{2k,1} \Sigma_{11}^{-1} x_{i(1)}$ are mutually uncorrelated.
We further assume they are independent. Intuitively, ${\bf z}_i$ can be viewed as
the unique part of irrelevant variables that cannot be explained by the relevant variables. Thus, for irrelevant variables, ${\bf z}_i$ is more representative than ${\bf x}_{i(2)}$. Assumption \ref{assump2} controls the lower bound of the nonzero parameters; see, for example, \citet{buhlmann2011}. Assumption \ref{assump3} imposes a lower bound, $N_1$, on the minimal eigenvalue of the covariance matrix of relevant variables. In practice, quantifying the rate under which $N_1$ decreases is difficult and problem specific, and it is frequently assumed to be constant; see, for example, \citet{Medeiros2016} and \citet{kock2015}. Assumption \ref{assump4} employs the strong irrepresentable condition of population covariance, which is similar to the condition in \citet{zhao2006}.

\indent To account for the cross-sectional dependence of the stationary process (${\bf x}_{i(1)}$) and (${\bf z}_i$), we also define the $\mathcal L^\infty$ functional dependence measure and its corresponding DAN, as follows:
\begin{eqnarray*}
&&\omega_{i,q,1}=\| \max_{1\leq j\leq s} |x_{ij}-x_{ij}^*| \|_q,\\
&&\| |{\bf x}_{.(1)}|_\infty \|_{q,\alpha}=\sup_{m\geq0}(m+1)^\alpha\Omega_{m,q,1},\quad \alpha\geq 0,\quad\text{and } \Omega_{m,q,1}=\sum_{i=m}^\infty \omega_{i,q,1}.
\end{eqnarray*}
Additionally, we define
\begin{eqnarray*}
\Psi_{q,\alpha,1}=\max_{1\leq j\leq s} \|x_{.j}\|_{q,\alpha}\quad\text{and }\quad \Upsilon_{q,\alpha,1} = \left(\sum_{j=1}^s \|x_{.j}\|_{q,\alpha}^q \right)^{1/q}.
\end{eqnarray*}
For $({\bf z}_i)$, the quantities $\| |{\bf z}_{.}|_\infty\|_{q,\alpha}$, $\Psi_{q,\alpha,2}$, and $\Upsilon_{q,\alpha,2}$ can be similarly defined. Clearly, $\Psi_{q,\alpha,1} \leq \| |{\bf x}_{.(1)}|_\infty\|_{q,\alpha}\leq \Upsilon_{q,\alpha,1}$ and $\Psi_{q,\alpha,2} \leq \| |{\bf z}_{.}|_\infty\|_{q,\alpha}\leq \Upsilon_{q,\alpha,2}$.

\indent Let $\sigma=\E e_i^2$. Define
\begin{eqnarray*}
\delta_*(\lambda,N_1,\sigma)&=&\frac{\lambda^2 s}{2n N_1}+\frac{2\sigma}{n},\\
M(\delta_*,\eta,\iota,\gamma)&=&\eta^{-1}\sqrt{\delta_* \log p}+\eta^{-1} n^{(\iota-1)/\gamma} \delta_*^{1/2}(\log p)^{3/2} \| |{\bf z}_{\cdot}|_\infty \|_{\gamma,\alpha_X}, \\
Q(\rho,\tau) &=& \sqrt{n\log s} + n^{\rho/\tau} (\log s)^{3/2} \| |{\bf x}_{\cdot(1)}|_\infty \|_{\gamma,\alpha_X}, \\
V_1(N_1) &=& \frac{s^2\log s}{N_1}, \\
V_2(N_1) &=& \frac{1}{N_1}s (\log s)^{3/2} \| |{\bf x}_{\cdot(1)}|_\infty \|_{\gamma,\alpha_X}^2 .
\end{eqnarray*}
These quantities are used in the following theorem.

\indent Theorem \ref{thm3} extends the results of \citet{zhao2006} to a random-design linear model with dependent errors. \citet{Medeiros2016} derived the asymptotic properties of sign consistency for the adaptive Lasso. In contrast, our results apply to the original Lasso, and do not need any assumptions on the weights. Note that even for heavy-tail variables, our results show that if the dependence among ${\bf z}_i$ is strong, the allowed dimension $p$ can be as large as some exponential of the sample size $n$; see Remark 2 for more details.

\begin{theorem} \label{thm3}
Suppose Assumptions \ref{assump1}, \ref{assump2}, \ref{assump3}, and \ref{assump4} hold. Assume that $\max_{1\leq j\leq p}\|x_{.j}\|_{\gamma,\alpha_X}<C_\gamma<\infty$ and $\|e_.\|_{q,\alpha_e}<C_q<\infty$, where $q,\gamma>4$, $\alpha_X, \alpha_e>0$, and constants $C_\gamma, C_q$ depend only on $\gamma,q$. Define
\[ \nu = \left\{\begin{array}{ll} 1 & \mbox{if $\alpha_X>1/2-2/\gamma$,} \\   \gamma/4-\alpha_X\gamma/2 & \mbox{if $\alpha_X <1/2-2/\gamma$,}
\end{array} \right. \]
and
\[ \iota = \left\{\begin{array}{ll} 1 & \mbox{if $\alpha_X>1/2-1/\gamma$,} \\   \gamma/2-\alpha_X\gamma & \mbox{if $\alpha_X <1/2-1/\gamma$.}
\end{array} \right. \]
 Let $\alpha=\min(\alpha_X,\alpha_e)$. Assume $\tau=q\gamma/(q+\gamma)>2$, and define
\[ \rho= \left\{\begin{array}{ll} 1 & \mbox{if $\alpha>1/2-1/\tau$,}  \\
\tau/2-\alpha\tau & \mbox{if $\alpha <1/2-1/\tau$.}\end{array}\right. \]
Furthermore, suppose $s=o(n)$. Then, for any $\lambda$ and sample size $n$, such that
\begin{eqnarray}
&&n \gtrsim V_1(N_1), \label{samplesize2} \\
&& n^{1-2\nu/\gamma} \gtrsim V_2(N_2), \label{samplesize3} \\
&&M(\delta_*,\eta,\iota,\gamma) + Q(\rho,\tau) \lesssim \lambda \le \frac{nN_1L}{4\sqrt{s}}, \label{regpar2}
\end{eqnarray}
the consistency probability $\P(\hat\beta =_s \beta)$ is at least
\begin{eqnarray}
1-C_1 (\log p)^{-\gamma}-C_2 (\log s)^{-\gamma/2}-C_3 (\log s)^{-\tau}-C_4 p^{-C_5}-C_6 s^{-C_7}-\frac{ \| e_{\cdot}\|_{q,\alpha_e}^q}{n^{q-1}\sigma^q} - \exp\left(-\frac{n\sigma^2}{\| e_{\cdot}\|_{2,\alpha_e}^2}\right). \label{signconsist}
\end{eqnarray}

\end{theorem}

\begin{rmk}
In particular, assume $N_1\asymp 1$, $\eta\asymp 1$. In addition, assume the weak temporal dependence case $\alpha_X>1/2-1/\gamma$ and $\alpha>1/2-1/\tau$. If the dependence measure $\| |{\bf x}_{.(1)}|_\infty\|_{\gamma,\alpha_X}\asymp s^{1/\gamma}$ and $\| |{\bf z}_{.}|_\infty\|_{\gamma,\alpha_X}\asymp p^{1/\gamma}$, which hold if all components $x_{ij}$ ($1\leq j\leq s$) and $z_{ik}$ ($s+1\leq k\leq p$) are nearly independent, then ($\ref{samplesize2}$), ($\ref{samplesize3}$), and ($\ref{regpar2}$) reduce to
\begin{eqnarray*}
n \gtrsim s^2\log s + s^{\frac{1+2/\gamma}{1-2/\gamma}}(\log s)^{\frac{3}{2-4/\gamma}} +s p^{2/\gamma} (\log p)^{3}
\end{eqnarray*}
and
\begin{eqnarray*}
\sqrt{n\log s}+n^{1/\tau}s^{1/\tau}(\log s)^{3/2} \lesssim \lambda \lesssim \frac{nL}{\sqrt{s}}.
\end{eqnarray*}
Additionally, if $s=O(n^{c_1})$, for some $c_1<\min\{1/2, (\gamma-2)/(\gamma+2)\}$, then the valid regularization parameter $\lambda$ has range $n^{1/2}+n^{1/\tau+c_1/\gamma} \ll \lambda \ll n^{1-c_1/2} L $. The dimension $p$ satisfies that $p \ll n^{\gamma(1-c_1)/2}$.

\indent On the other hand, assume $\| |{\bf x}_{.(1)}|_\infty\|_{\gamma,\alpha_X} \asymp s^{1/\gamma}$ and $\| |{\bf z}_{.}|_\infty\|_{\gamma,\alpha_X} \asymp 1$; that is, all components $z_{ik}$ ($s+1\leq k\leq p$) are strongly dependent. Let $s=O(n^{c_1})$, for some $c_1<\min\{1/2, (\gamma-2)/(\gamma+2)\}$. Then the existence of regularization parameter $\lambda$ requires $n^{1/2}+n^{1/\tau+c_1/\gamma} \ll \lambda \ll n^{1-c_1/2} L $. The dimension $p$ satisfies $p \ll \exp\{n^{(1-c_1)/3}\}$.

\indent Furthermore, if $\| |{\bf x}_{.(1)}|_\infty\|_{\gamma,\alpha_X} \asymp 1$ and $\| |{\bf z}_{.}|_\infty\|_{\gamma,\alpha_X} \asymp 1$, and $s=O(n^{c_1})$ for some $c_1< 1/2$, then the existence of regularization parameter $\lambda$
requires $n^{1/2} \ll \lambda \ll n^{1-c_1/2} L $, and the dimension $p$
satisfies  $p \ll \exp\{n^{(1-c_1)/3}\}$.

\indent In summary, the allowed dimension $p$ varies from $n^{\gamma(1-c_1)/2}$ to $\exp\{n^{(1-c_1)/3}\}$, depending on the cross-sectional dependence of $z_{ik}, s+1\le k\le p$. 
\end{rmk}

\indent Note that if the assumptions in Example \ref{example1} hold, then together with the strong irrepresentable condition, the results of Theorem \ref{thm3} continue to apply. In general, the strong irrepresentable condition is nontrivial, particularly because we do not know sign$(\beta)$ a priori. Then, we need the strong irrepresentable condition to hold for every possible combination of signs and placement of zeros. We give a simple example below in which the strong irrepresentable condition is guaranteed. All diagonal elements of $\Sigma$ are assumed to be one which is equivalent to normalizing all covariates in the model to the same scale, because the strong irrepresentable condition is invariant under any common scaling of $\Sigma$.

\begin{exmp}
Consider the following autoregressive model with exogenous variables:
\begin{eqnarray}
y_i=\sum_{l=1}^a\phi_{l}y_{i-l}+\psi{\bf z}_{i}+e_i=\beta'{\bf x}_i+e_i,
\end{eqnarray}
where $a$ is nonnegative finite integer, ${\bf z}_i$ is independent of $e_i$, and the errors $e_i$ are homogeneous. Assume the roots of the polynomial $1-\sum_{l=1}^a\phi_l B^l$ are outside the unit circle, which ensures the stationarity of the autoregressive part of the model. In addition, assume $\Sigma=\E{\bf x}_i{\bf x}_i'$ is positive definite.\\
\indent Furthermore, suppose $\beta$ has $s$ nonzero entries. Similarly to Corollary 2 in \citet{zhao2006}, if $\Sigma$ has ones on the diagonal and the bounded correlation $|\sigma_{jk}|\leq c/(2s-1)$, for a constant $0< c< 1$, then the strong irrepresentable condition holds. In this case, we need the autocorrelation of $y_i$ to be weak,
and all covariates ${\bf z}_i$ are slightly correlated.
\end{exmp}

\begin{rmk}
The Lasso may fail in the presence of strong serial dependence.
Consider two scalar Gaussian autoregressive,
AR(3), models:
\begin{eqnarray}\label{armodel1}
y_i=1.9 y_{i-1}-0.8 y_{i-2}-0.1 y_{i-3}+e_i,
\end{eqnarray}
and
\begin{eqnarray}\label{armodel2}
y_i= y_{i-1}-0.8 y_{i-2}-0.1 y_{i-3}+e_i,
\end{eqnarray}
where $e_i$ follows the standard normal distribution. Then, AR(3) in model (\ref{armodel1}) is unit-root nonstationary, but that in model (\ref{armodel2}) is stationary. We generate 2000 observations from each of the two models. We choose $y_{i-10},y_{i-9},...,y_{i-1}$ and $x_{1i},...,x_{10,i}$ as regressors, where $x_{li}$ are i.i.d. standard normal.
Figure~\ref{fig:lasso} shows the model selection results for scaling versus not scaling the predictors.
\begin{figure}
\begin{center}
\subfloat[scaling for AR model (\ref{armodel1})]{\includegraphics[scale=0.6]{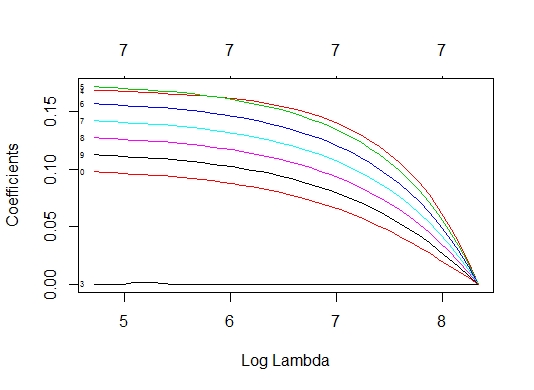}}
\subfloat[not scaling for AR model (\ref{armodel1})]{\includegraphics[scale=0.6]{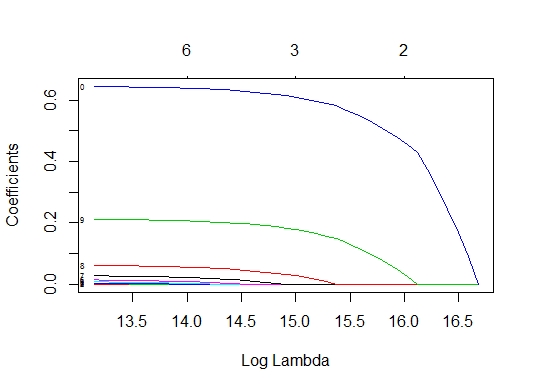}}\\
\subfloat[scaling for AR model (\ref{armodel2})]{\includegraphics[scale=0.6]{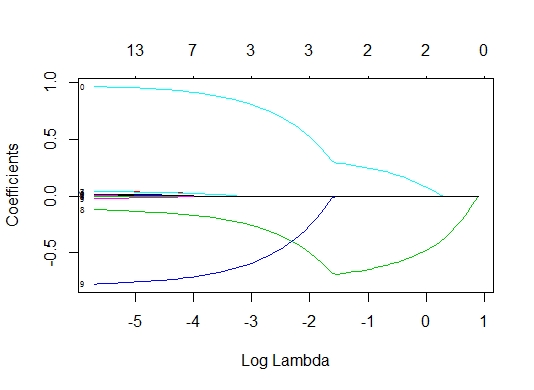}}
\subfloat[not scaling for AR model (\ref{armodel2})]{\includegraphics[scale=0.6]{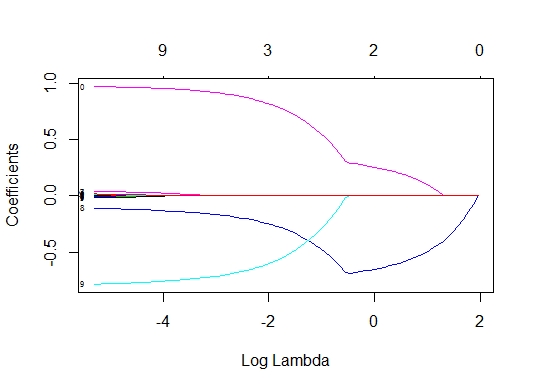}}
\end{center}
\caption{Results of Lasso regression for the two AR(3) series in (\ref{armodel1}) and
(\ref{armodel2}), created using the glmnet package of R.}
\label{fig:lasso}
\end{figure}

The default Lasso procedure standardizes each variable in $y_i$. For unit-root nonstationary time series, standardization might
wash out the dependence of the stationary part; see (a) and (b) of
Figure~\ref{fig:lasso}. In this paper, we only consider stationary time series for which scaling the predictors
does not affect the estimation consistency of the Lasso estimates; see
(c) and (d) of Figure~\ref{fig:lasso}.

The following proposition shows a necessary and sufficient condition for a stationary AR(2) model, under which the strong irrepresentable condition (Assumption \ref{assump4}) holds. Similar results hold for the general stationary AR($d$) model.
\begin{proposition}
Consider the stationary AR(2) model,
$$y_i=\phi_1 y_{i-1} + \phi_2 y_{i-2}+ e_i,$$
where $e_i$ are i.i.d. random variates with mean zero and finite variance.
We also normalize $y_i$, such that the variance of $y_i$ is one.
Then, the strong irrepresentable condition (Assumption \ref{assump4}) holds if and only if
\begin{eqnarray}
|\phi_1|+|\phi_2|<1.
\end{eqnarray}

\end{proposition}

\end{rmk}

\section{Simulation Study}
In this section, we use a simulation to demonstrate the performance of the Lasso regression for
dependent data in finite samples, and to compare its efficacy with that of the mixed-frequency data
sampling regression (MIDAS) commonly used in the econometric literature; see \citet{ghysels2004}. In addition,
we compare the model selection consistency and parameter estimation of the Lasso estimator and the Dantzig estimator for dependent data in finite samples.

We first consider the following data-generating process:
\begin{eqnarray}\label{simu}
&& y_i=\phi y_{i-1}+\bm x_{i-1,1}^T\beta_s+e_i, \notag\\
&& \bm x_i=\begin{bmatrix} \bm x_{i,1}\\ \bm x_{i,2} \end{bmatrix}=\sum_{j=1}^m A_j \begin{bmatrix} \bm x_{i-j,1}\\ \bm x_{i-j,2} \end{bmatrix} + \bm \eta_i,
\end{eqnarray}
where $\phi=0.6$, each element of $\beta_s$ is given by $\beta_{s,j}=\frac{1}{\sqrt s}(-1)^{j}$, and $\bm x_{i,1}$ is an $s\times 1$ vector of relevant variables. Let $\beta=(\beta_s,\beta_{s^c})$, where $\beta_{s^c}=\bm 0$ is a $(p-s)\times 1$ vector. The errors $e_i$ and $\bm \eta_{ij}$ are i.i.d. random variables from a Student-$t$ distribution with five degrees of
freedom,  and $e_i$ and $\bm \eta_i$ are all mutually uncorrelated. The explanatory variable  process $\bm x_i$, which has $p-s$ irrelevant variables, follows a vector autoregressive, VAR($m$), model. The following two choices of $\bm x_i$ are considered, denoted as Model 1 and Model 2, respectively.
\begin{enumerate}[(1).]
\item {\bf Model 1}: The explanatory process $\bm x_i$ is a VAR(4) process, where $A_1$ and $A_4$ assume a block-diagonal structure, and $A_2=A_3=0$. In particular, the first two and the last two blocks are $5\times 5$ matrices, with all entries of the blocks of $A_1$ equal to $0.15$, and all entries of the blocks of $A_4$ equal to $-0.1$. The other blocks are $10\times 10$ matrices, with all elements of the blocks of $A_1$ equal to $0.075$, and all elements of the blocks of $A_4$ equal to $-0.05$. This structure could be motivated by a model built for mixed-frequency data with some quarterly time series, often encountered in macroeconomic analysis.
\item {\bf Model 2}: The explanatory process $\bm x_i$ follows a VAR(1) model, where $A_1$ is block-diagonal, with the same block structure given by Model 1. The $(j,k)$th entry of the block is $(-1)^{|j-k|}\rho^{|j-k|+1}$, with $\rho=0.4$. Hence, the entries decrease exponentially fast with their distances from the diagonal.
\end{enumerate}

\indent We employ sample sizes $n=50, 100, 200$, with different choices of $p$ and $s$. We set $p=100, 200, 400$ and $s=5, 10, 20$. For comparison, we also simulate a response series from a MIDAS model. In Model (\ref{simu}), for $s=5, 10, 20$, let $\beta_s =\beta(1)$, $(\beta(1)^T,\beta(2)^T)^T$ or $(\beta(1)^T,\beta(2)^T,\beta(3)^T)^T$, respectively, with
\begin{eqnarray}\label{weight}
\beta_j(l)=\frac{\exp(\delta_1j+\delta_2 j^2)}{\sum_{k=1}^{|\beta(l)|_0}\exp(2\delta_1 k+2\delta_2 k^2)},
\end{eqnarray}
where $\beta(1)$ and $\beta(2)$ have five variables, $\beta(3)$ has 10 variables, and $\delta=(\delta_1,\delta_2)'=(0.5,-1)'$. All other settings remain the same.
The two choices of $\bm x_i$ in Models 1 and 2 are used, and we denote the
resulting MIDAS models as Models 3 and 4, respectively.
 The models estimated by the Lasso have $\lambda$ selected using the BIC; see \citet{buhlmann2011}. The consistency of the Lasso estimator selected by the BIC was first proved by \citet{zou2007} under the case $p<n$. Then, \citet{tibshirani2012} studied the effective degrees of freedom of the Lasso when $p>n$. It would be interesting to investigate the theoretical justification of the consistency of the BIC criterion for the Lasso under the time series setting. We leave this to future work. We also employed models with $\lambda$ selected using cross-validation. However, cross-validation does not improve the results, and is considerably slower in its computation. For the models estimated by MIDAS, we only consider the exponential Almon lag polynomial weighting scheme (see (\ref{weight})) for the first 100 variables, and impute the true values as initial values.

\indent Table \ref{tab:table1} shows the average absolute error (AE) and average root mean squared error (RMSE) for the Lasso estimators and MIDAS estimators over the 10,000 Monte Carlo simulations for the data-generating processes used. The AE and the RMSE are defined as
\begin{eqnarray*}
\text{AE}&=& \frac{1}{MC}\sum_{l=1}^{MC} |(\hat\phi;\hat\beta)-(\phi;\beta)|_1 , \\
\text{RMSE}&=&\sqrt{\frac{1}{MC}\sum_{l=1}^{MC} |(\hat\phi;\hat\beta)-(\phi;\beta)|_2^2 },
\end{eqnarray*}
where MC denotes the number of Monte Carlo repetitions.
From the table, it is clear that both measures show that the Lasso regression
provides a substantially more accurate parameter estimation than that of the mixed-frequency data sampling regression (MIDAS) in the presence of irrelevant variables. Furthermore, as expected, the AE and the RMSE of the estimators decrease with $n$,  but increase with $s$ and $p$.

\indent To evaluate the performance of out-of-sample forecasts, we use the estimated parameters to compute one-step-ahead forecasts, and consider 10 out-of-sample predictions, denoted by $y_{n+1}, \ldots, y_{n+10}$. Table \ref{tab:table2} shows the average absolute forecast error (AFE) and average root mean squared forecast error (RMSFE) over the 10,000 Monte Carlo simulations, which are calculated as
\begin{eqnarray*}
\text{AFE}&=& \frac{1}{10MC}\sum_{l=1}^{MC} \sum_{k=1}^{10} |\hat y_{n+k}-y_{n+k}| , \\
\text{RMSFE}&=&\sqrt{ \frac{1}{10MC}\sum_{l=1}^{MC} \sum_{k=1}^{10} |\hat y_{n+k}-y_{n+k}|^2}.
\end{eqnarray*}
The forecasting results in Table \ref{tab:table2} show that the Lasso regression has smaller AE and RMSFE than those of the MIDAS in all settings. Furthermore, the results show clearly that the performance of the Lasso regression and the MIDAS improves with the sample size, but deteriorates as the number of relevant variables $s$ increases. Finally, both the AE and the RMSFE of the Lasso regression decrease faster than those of the MIDAS as the sample size $n$ increases.
In fact, the AE and RMSFE of the MIDAS remain high even when $n=200$. Because we only fit the MIDAS through the first 100 variables, its performance does not change as $p$ increases.
Overall, in the presence of irrelevant variables, the Lasso regression significantly outperforms the MIDAS regression.

\begin{table}[htbp]
\renewcommand{\arraystretch}{0.59}
\begin{widepage}
\centering
\caption{Accuracy in parameter estimation of Lasso regression and mixed-frequency data sampling
regression. The results are based on 10,000 repetitions, where AE and RMSE denote the
average mean absolute errors and average root mean squared errors over Monte Carlo
repetitions and parameters. In the table, $s$, $p$, and $n$ denote the number of nonzero
parameters, dimension of regressors, and sample size, respectively. }
\label{tab:table1}
\scriptsize
\begin{tabular}{cccccccccccccccccc}\hline
$s$ & $n$ & & \multicolumn{7}{c}{Absolute Error (AE) $\times 10^{2}$} & & \multicolumn{7}{c}{Root Mean Square Error (RMSE) $\times 10^{2}$}\\ \cline{4-10}\cline{12-18}
  & & &\multicolumn{3}{c}{Lasso} & & \multicolumn{3}{c}{MIDAS} & & \multicolumn{3}{c}{Lasso} & & \multicolumn{3}{c}{MIDAS}\\ \cline{4-6}\cline{8-10}\cline{12-14}\cline{16-18}
  & & & \multicolumn{15}{c}{$p$}\\ \cline{4-18}
  & & & 100 & 200 & 400 & & 100 & 200 & 400 & & 100 & 200 & 400 & & 100 & 200 & 400\\ \hline
  \multicolumn{5}{l}{Model 1}\\ \hline
  \multirow{3}{*}{5}  & 50 & &  2.44 & 2.64 & 2.84 & & 6.63 & 6.64 & 6.67 & & 3.08 & 3.75 & 4.55 & & 3.73 & 4.44 & 5.29\\
  & 100 & & 1.89 & 2.07 & 2.22 & & 6.24 & 6.26 & 6.28 & & 2.79 & 3.44 & 4.21 & & 3.64 & 4.33 & 5.15\\
  & 200 & & 1.27 & 1.49 & 1.70 & & 5.91 & 5.91 & 5.94 & & 2.28 & 2.92 & 3.71 & & 3.56 & 4.23 & 5.04\\
  \multirow{3}{*}{10}  & 50 & & 4.60 & 4.99 & 5.30 & & 8.26 & 8.31 & 8.32 & & 3.66 & 4.45 & 5.38 & & 4.09 & 4.88 & 5.80\\
  & 100 & & 3.69 & 4.11 & 4.39 & & 7.86 & 7.88 & 7.90 & & 3.36 & 4.17 & 5.10 & & 4.02 & 4.78 & 5.69\\
  & 200 & & 2.28 & 2.74 & 3.29 & & 7.50 & 7.55 & 7.56 & & 2.65 & 3.39 & 4.42 & & 3.96 & 4.71 & 5.60\\
  \multirow{3}{*}{20}  & 50 & & 7.83 & 8.81 & 8.93 & & 10.76 & 10.82 & 10.83 & & 4.08 & 5.00 & 6.00 & & 4.42 & 5.26 & 6.26\\
  & 100 & & 6.56 & 7.33 & 7.70 & & 10.38 & 10.43 & 10.44 & & 3.84 & 4.75 & 5.77 & & 4.35 & 5.18 & 6.16\\
  & 200 & & 4.69 & 5.55 & 6.56 & & 10.08 & 10.12 & 10.15 & & 3.31 & 4.21 & 5.40 & & 4.30 & 5.12 & 6.09\\ \hline
  \multicolumn{5}{l}{Model 2}\\ \hline
  \multirow{3}{*}{5}  & 50 & & 0.95 & 1.14 & 1.38 & & 4.95 & 4.97 & 4.99 & & 2.02 & 2.53 & 3.19 & & 3.31 & 3.94 & 4.70\\
  & 100 & & 0.54 & 0.60 & 0.67 & & 4.55 & 4.58 & 4.58 & & 1.56 & 1.92 & 2.36 & & 3.18 & 3.79 & 4.50\\
  & 200 & & 0.34 & 0.36 & 0.38 & & 4.20 & 4.21 & 4.22 & & 1.26 & 1.53 & 1.87 & & 3.06 & 3.64 & 4.33\\
  \multirow{3}{*}{10}  & 50 & & 1.91 & 2.40 & 2.92 & & 5.46 & 5.46 & 5.46 & & 2.54 & 3.26 & 4.18 & & 3.53 & 4.20 & 4.99\\
  & 100 & & 1.06 & 1.24 & 1.46 & & 5.03 & 5.07 & 5.08 & & 1.92 & 2.41 & 3.04 & & 3.39 & 4.04 & 4.81\\
  & 200 & & 0.65 & 0.71 & 0.79 & & 4.60 & 4.63 & 4.65 & & 1.52 & 1.87 & 2.31 & & 3.25 & 3.87 & 4.61\\
  \multirow{3}{*}{20}  & 50 & & 3.19 & 4.21 & 4.94 & & 6.12 & 6.15 & 6.18 & & 2.95 & 3.85 & 4.96 & & 3.76 & 4.48 & 5.34\\
  & 100 & & 1.75 & 2.14 & 2.59 & & 5.68 & 5.69 & 5.70 & & 2.23 & 2.85 & 3.64 & & 3.63 & 4.32 & 5.15\\
  & 200 & & 1.07 & 1.21 & 1.38 & & 5.26 & 5.27 & 5.29 & & 1.77 & 2.20 & 2.74 & & 3.51 & 4.18 & 4.98\\ \hline
  \multicolumn{5}{l}{Model 3}\\ \hline
  \multirow{3}{*}{5}  & 50 & & 1.71 & 2.05 & 2.43 & & 6.57 & 6.60 & 6.62 & & 2.62 & 3.29 & 4.11 & & 3.74 & 4.46 & 5.30\\
  & 100 & & 0.93 & 1.06 & 1.21 & & 6.27 & 6.31 & 6.33 & & 2.03 & 2.54 & 3.18 & & 3.65 & 4.35 & 5.18\\
  & 200 & & 0.57 & 0.63 & 0.69 & & 6.17 & 6.20 & 6.21 & & 1.62 & 2.02 & 2.50 & & 3.61 & 4.30 & 5.11\\
  \multirow{3}{*}{10}  & 50 & & 3.74 & 4.47 & 5.07 & & 8.41 & 8.44 & 8.46 & & 3.34 & 4.17 & 5.16 & & 4.13 & 4.92 & 5.86\\
  & 100 & & 2.06 & 2.52 & 3.00 & & 8.20 & 8.24 & 8.25 & & 2.59 & 3.32 & 4.25 & & 4.08 & 4.85 & 5.77\\
  & 200 & & 1.20 & 1.38 & 1.58 & & 8.10 & 8.14 & 8.16 & & 2.00 & 2.52 & 3.18 & & 4.05 & 4.82 & 5.73\\
  \multirow{3}{*}{20}  & 50 & & 7.23 & 8.77 & 9.38 & & 11.02 & 11.07 & 11.09 & & 3.90 & 4.88 & 5.95 & & 4.47 & 5.32 & 6.32\\
  & 100 & & 4.45 & 5.81 & 7.01 & & 10.92 & 10.97 & 11.00 & & 3.22 & 4.16 & 5.32 & & 4.43 & 5.28 & 6.28\\
  & 200 & & 2.53 & 2.93 & 3.50 & & 10.87 & 10.93 & 10.95 & & 2.49 & 3.11 & 3.97 & & 4.42 & 5.26 & 6.26\\ \hline
  \multicolumn{5}{l}{Model 4}\\ \hline
  \multirow{3}{*}{5}  & 50 & & 1.39 & 1.58 & 1.78 & & 5.14 & 5.16 & 5.16 & & 2.49 & 3.10 & 3.83 & & 3.47 & 4.13 & 4.90\\
  & 100 & & 0.96 & 1.05 & 1.12 & & 4.58 & 4.59 & 4.59 & & 2.12 & 2.63 & 3.22 & & 3.31 & 3.95 & 4.69\\
  & 200 & & 0.71 & 0.77 & 0.83 & & 4.22 & 4.23 & 4.25 & & 1.83 & 2.28 & 2.81 & & 3.23 & 3.85 & 4.58\\
  \multirow{3}{*}{10}  & 50 & & 2.40 & 2.79 & 3.14 & & 6.03 & 6.07 & 6.10 & & 2.90 & 3.64 & 4.54 & & 3.80 & 4.53 & 5.39\\
  & 100 & & 1.67 & 1.86 & 2.02 & & 5.50 & 5.53 & 5.55 & & 2.47 & 3.08 & 3.80 & & 3.69 & 4.39 & 5.23\\
  & 200 & & 1.23 & 1.38 & 1.50 & & 5.13 & 5.15 & 5.16 & & 2.11 & 2.65 & 3.30 & & 3.62 & 4.31 & 5.13\\
  \multirow{3}{*}{20}  & 50 & & 3.68 & 4.58 & 4.97 & & 6.88 & 6.93 & 6.93 & & 3.22 & 4.09 & 5.14 & & 4.06 & 4.84 & 5.75\\
  & 100 & & 2.43 & 2.77 & 3.06 & & 6.38 & 6.42 & 6.45 & & 2.71 & 3.41 & 4.25 & & 3.96 & 4.72 & 5.62\\
  & 200 & & 1.78 & 2.00 & 2.22 & & 6.04 & 6.08 & 6.08 & & 2.32 & 2.91 & 3.66 & & 3.91 & 4.66 & 5.55\\ \hline
\end{tabular}
\end{widepage}
\end{table}

\begin{table}[htbp]
\renewcommand{\arraystretch}{0.59}
\begin{widepage}
\centering
\caption{Performance of Out-of-sample predictions of Lasso regression and mixed frequency
data sampling regression (MIDAS). The results are based on 10 one-step ahead predictions and
10,000 iterations, where AFE and RMSFE denote the average absolute forecast errors and
root mean squared forecast errors, respectively, and $s$, $p$, and $n$ are the
number of nonzero parameters, dimension of regressors, and sample size, respectively. For MIDAS,
the maximum $p$ is fixed at 100.}
\label{tab:table2}
\scriptsize
\begin{tabular}{ccrrrrrrrrrrrrrrrr} \hline 
$s$ & $n$ & & \multicolumn{7}{c}{Absolute Error (AE) $\times 10^{2}$} & & \multicolumn{7}{c}{Root Mean Square Forecast Error (RMSFE) $\times 10^{2}$}\\ \cline{4-10}\cline{12-18}
  & & &\multicolumn{3}{c}{Lasso} & & \multicolumn{3}{c}{MIDAS} & & \multicolumn{3}{c}{Lasso} & & \multicolumn{3}{c}{MIDAS}\\ \cline{4-6}\cline{8-10}\cline{12-14}\cline{16-18}
  & & & \multicolumn{15}{c}{$p$}\\ \cline{4-18}
  & & & 100 & 200 & 400 & & 100 & 200 & 400 & & 100 & 200 & 400 & & 100 & 200 & 400\\ \hline
  \multicolumn{5}{l}{Model 1}\\ \hline
\multirow{3}{*}{5}  & 50 & & 120.0 & 125.8 & 130.3 & & 169.2 & 161.5 & 162.3 & & 147.2 & 153.8 & 158.8 & & 206.0 & 197.4 & 198.0\\
  & 100 & & 102.7 & 106.6 & 110.7 & & 162.2 & 156.7 & 156.4 & & 127.7 & 132.3 & 136.9 & & 197.7 & 191.7 & 191.3 \\
  & 200 & & 86.9 & 90.6 & 95.4 & & 156.8 & 152.4 & 153.2 & & 109.7 & 114.1 & 119.4 & & 191.4 & 186.7 & 187.4 \\
  \multirow{3}{*}{10}  & 50 & & 151.6 & 159.6 & 166.4 & & 185.0 & 178.8 & 179.9 & & 185.2 & 194.3 & 202.0 & & 225.9 & 218.6 & 219.9\\
  & 100 & & 125.6 & 133.9 & 141.7 & & 177.6 & 171.6 & 172.3 & & 155.1 & 164.9 & 173.8 & & 216.9 & 210.1 & 211.4 \\
  & 200 & & 96.0 & 101.9 & 112.2 & & 171.9 & 167.7 & 168.2 & & 120.3 & 127.2 & 139.3 & & 210.2 & 206.1 & 206.3 \\
  \multirow{3}{*}{20}  & 50 & & 177.7 & 188.8 & 195.0 & & 205.0 & 200.0 & 199.9 & & 216.1 & 229.2 & 236.2 & & 250.2 & 244.5 & 244.1\\
  & 100 & & 150.2 & 162.2 & 170.0 & & 195.8 & 191.7 & 191.1 & & 184.0 & 198.5 & 207.7 & & 239.5 & 235.0 & 234.4 \\
  & 200 & & 118.6 & 128.7 & 145.1 & & 190.1 & 185.9 & 188.2 & & 146.7 & 159.2 & 178.4 & & 232.4 & 228.4 & 230.5 \\ \hline
  \multicolumn{5}{l}{Model 2}\\ \hline
  \multirow{3}{*}{5} & 50 & & 96.4 & 101.8 & 107.2 & & 147.3 & 148.8 & 148.5 & & 119.7 & 125.5 & 131.7 & & 179.9 & 181.5 & 180.9\\
  & 100 & & 84.1 & 85.7 & 88.1 & & 142.1 & 142.7 & 142.9 & & 106.2 & 108.1 & 110.4 & & 173.3 & 174.0 & 174.0 \\
  & 200 & & 78.4 & 79.6 & 80.6 & & 138.6 & 137.6 & 138.8 & & 99.9 & 101.5 & 102.3 & & 169.0 & 168.1 & 169.2 \\
  \multirow{3}{*}{10}  & 50 & & 114.1 & 125.7 & 140.0 & & 171.5 & 164.2 & 163.9 & & 139.9 & 153.7 & 169.8 & & 208.0 & 199.6 & 199.5\\
  & 100 & & 90.9 & 95.3 & 100.8 & & 156.7 & 157.9 & 158.0 & & 114.1 & 118.7 & 124.9 & & 190.6 & 191.9 & 191.9 \\
  & 200 & & 81.7 & 83.1 & 85.3 & & 151.4 & 151.1 & 151.9 & & 103.7 & 105.1 & 107.6 & & 184.1 & 183.7 & 184.5 \\
  \multirow{3}{*}{20}  & 50 & & 126.9 & 144.5 & 167.8 & & 178.2 & 173.1 & 173.5 & & 155.4 & 175.9 & 202.9 & & 216.5 & 211.1 & 211.6\\
  & 100 & & 97.7 & 105.1 & 113.7 & & 169.7 & 164.3 & 164.7 & & 121.6 & 130.1 & 139.9 & & 206.5 & 200.4 & 200.9 \\
  & 200 & & 85.3 & 87.9 & 91.9 & & 161.7 & 157.1 & 158.0 & & 107.8 & 110.5 & 115.1 & & 196.9 & 191.9 & 192.9 \\ \hline
  \multicolumn{5}{l}{Model 3}\\ \hline
  \multirow{3}{*}{5}  & 50 & & 117.4 & 128.7 & 140.5 & & 152.9 & 153.1 & 153.3 & & 143.0 & 155.5 & 168.3 & & 187.5 & 187.6 & 187.9\\
  & 100 & & 89.4 & 92.8 & 97.1 & & 144.7 & 145.0 & 145.0 & & 112.2 & 116.0 & 120.4 & & 144.7 & 178.2 & 178.1 \\
  & 200 & & 80.6 & 81.2 & 82.8 & & 142.3 & 141.0 & 141.1 & & 102.3 & 103.0 & 105.0 & & 174.7 & 173.6 & 173.5 \\
  \multirow{3}{*}{10}  & 50 & & 154.4 & 172.5 & 188.9 & & 178.9 & 179.1 & 179.8 & & 185.9 & 206.0 & 224.6 & & 218.4 & 218.8 & 219.7\\
  & 100 & & 103.1 & 112.9 & 124.3 & & 171.2 & 171.3 & 170.6 & & 127.9 & 138.7 & 152.2 & & 209.7 & 209.5 & 209.1 \\
  & 200 & & 84.7 & 88.2 & 91.0 & & 166.9 & 168.2 & 167.4 & & 107.1 & 111.1 & 114.3 & & 204.7 & 205.9 & 205.3 \\
  \multirow{3}{*}{20}  & 50 & & 197.3 & 224.5 & 244.3 & & 206.9 & 205.0 & 205.3 & & 236.1 & 266.5 & 288.4 & & 251.6 & 249.7 & 249.8\\
  & 100 & & 130.9 & 150.8 & 172.2 & & 196.6 & 197.6 & 196.4 & & 160.2 & 182.9 & 207.6 & & 240.2 & 240.9 & 239.3 \\
  & 200 & & 97.0 & 101.2 & 109.0 & & 193.1 & 193.8 & 193.7 & & 121.2 & 125.9 & 134.8 & & 236.5 & 237.4 & 237.1 \\ \hline
  \multicolumn{5}{l}{Model 4}\\ \hline
  \multirow{3}{*}{5}  & 50 & & 103.0 & 108.7 & 113.2 & & 131.7 & 131.9 & 130.9 & & 126.8 & 133.3 & 138.2 & & 162.6 & 162.9 & 161.7\\
  & 100 & & 88.4 & 90.4 & 92.9 & & 121.6 & 122.0 & 121.5 & & 110.9 & 113.0 & 115.8 & & 150.9 & 151.5 & 150.8 \\
  & 200 & & 81.3 & 82.6 & 83.4 & & 118.0 & 117.1 & 116.8 & & 103.3 & 104.4 & 105.4 & & 147.0 & 145.7 & 145.2 \\
  \multirow{3}{*}{10}  & 50 & & 117.6 & 126.6 & 136.2 & & 148.6 & 148.5 & 148.5 & & 144.1 & 154.4 & 165.7 & & 183.7 & 183.0 & 183.0\\
  & 100 & & 95.8 & 99.8 & 103.4 & & 139.4 & 139.5 & 139.3 & & 119.8 & 124.2 & 128.2 & & 172.5 & 172.6 & 172.3 \\
  & 200 & & 84.9 & 87.5 & 89.7 & & 134.2 & 135.0 & 134.7 & & 107.3 & 110.1 & 112.7 & & 166.5 & 167.3 & 167.0 \\
  \multirow{3}{*}{20}  & 50 & & 132.2 & 148.5 & 162.3 & & 163.8 & 164.7 & 163.7 & & 161.3 & 180.2 & 196.3 & & 201.7 & 202.7 & 201.6\\
  & 100 & & 102.4 & 108.9 & 115.4 & & 154.2 & 154.0 & 154.7 & & 127.2 & 134.8 & 142.3 & & 190.6 & 190.7 & 190.9 \\
  & 200 & & 88.6 & 92.1 & 96.2 & & 150.0 & 149.8 & 150.2 & & 111.7 & 115.7 & 120.2 & & 185.8 & 185.4 & 185.7 \\ \hline
\end{tabular}
\end{widepage}
\end{table}

Next, we compare the model selection and parameter estimation of the Lasso estimator and the Dantzig estimator for dependent data. We use the same data-generating process (\ref{simu}), where $\phi=0.6$ and $\bm x_{i,1}$ is an $s\times 1$ vector of relevant variables. Here, we set each element of $\beta_s$  by $\beta_{s,j}=3(-1)^{j}$. Model 1 and Model 2, defined previously, are chosen for $\bm x_i$. Table \ref{tab:tablec} shows the number of noise covariates that are selected (False Positive), number of signal covariates that are not
selected (False Negative), and average root mean squared error (RMSE) for the Lasso estimators and the Dantzig estimators over the 10,000 Monte Carlo simulations for the data-generating processes used. As expected, False Positive and RMSE decrease with $n$,  but increase with $s$ and $p$. False Negative for the two methods are almost the same. In terms of False Negative and RMSE, the Lasso estimator substantially outperforms the
Dantzig selector. The Dantzig selector might be more sensitive to heavy tails and outliers,
because it uses the $L_\infty$ norm.
The rate of convergence of the Lasso estimator in our study is faster than that of the Dantzig selector in \citet{wu2016}. They built an $L_\infty$-type rate of convergence for the Dantzig estimator, which is related to the unknown $L_1$ norm of the true coefficients and the matrix $L_1$ norm of the population matrix. We overcome this weakness and achieve the same bounds for the Lasso regression under i.i.d. data, but with different requirements for the regularization parameter $\lambda$ and sample size $n$.

\begin{table}[htbp]
\renewcommand{\arraystretch}{0.59}
\begin{widepage}
\centering
\caption{Accuracy in model selection and parameter estimation of Lasso estimator and Dantzig estimator for linear regression. The results are based on 10,000 repetitions, where RMSE denote the average root mean squared error over Monte Carlo repetitions and parameters. In the table, $s$, $ p$, and $n$ denote the number of nonzero parameters, dimension of regressors, and sample size, respectively. }
\label{tab:tablec}
\scriptsize
\begin{tabular}{cccccccccccccccccc}\hline
$s$ & $n$ & & \multicolumn{7}{c}{Model 1} & & \multicolumn{7}{c}{Model 2}\\ \cline{4-10}\cline{12-18}
  & & &\multicolumn{3}{c}{Lasso} & & \multicolumn{3}{c}{Dantzig} & & \multicolumn{3}{c}{Lasso} & & \multicolumn{3}{c}{Dantzig}\\ \cline{4-6}\cline{8-10}\cline{12-14}\cline{16-18}
  & & & \multicolumn{15}{c}{$p$}\\ \cline{4-18}
  & & & 100 & 200 & 400 & & 100 & 200 & 400 & & 100 & 200 & 400 & & 100 & 200 & 400\\ \hline
  \multicolumn{5}{l}{False Negative}\\ \hline
  \multirow{3}{*}{5}  & 50 & &  0.077 & 0.20 & 0.67 & & 0.072 & 0.28 & 0.73 & & 0 & 0 & 0.003 & & 0 & 0 & 0.002 \\
  & 100 & & 0 & 0 & 0.01 & & 0 & 0  & 0.04  & & 0 & 0 & 0 & & 0 & 0 & 0 \\
  & 200 & & 0 & 0 & 0 & & 0 & 0  & 0  & & 0 & 0 & 0 & & 0 & 0 & 0 \\
  \multirow{3}{*}{10}  & 50 & & 0.67 & 2.07 & 4.28 & & 0.81 & 2.50 & 4.04 & & 0.011 & 0.14 & 0.79 & & 0.045 & 0.17 & 0.95 \\
  & 100 & & 0 & 0.004 & 0.11 & & 0 & 0.006 & 0.15  & & 0 & 0 & 0 & & 0 & 0 & 0 \\
  & 200 & & 0 & 0 & 0 & & 0 & 0 & 0  & & 0 & 0 & 0 & & 0 & 0 & 0 \\
  \multirow{3}{*}{20}  & 50 & & 4.65 & 7.45 & 9.61 & & 4.83 & 7.29 & 10.1 & & 1.94 & 5.70 & 7.18 & & 2.48 & 5.45 & 8.34 \\
  & 100 & & 0 & 0.21 & 2.27 & & 0.03 & 0.28 & 2.14 & & 0 & 0 & 0.029 & & 0 & 0.002 & 0.052 \\
  & 200 & & 0 & 0 & 0 & & 0 & 0 & 0 & & 0 & 0 & 0 & & 0 & 0 & 0 \\ \hline
  \multicolumn{5}{l}{False Positive}\\ \hline
  \multirow{3}{*}{5}  & 50 & & 10.9 & 15.0 & 22.5 & & 15.85 & 24.20 & 32.0 & & 5.25 & 8.01 & 11.5 & & 6.45 & 13.3 & 21.7 \\
  & 100 & & 5.30 & 8.97 & 13.6 & & 8.03 & 14.11 & 19.8 & & 2.05 & 3.02 & 4.00 & & 4.14 & 7.43 & 9.91\\
  & 200 & & 1.63 & 3.01 & 4.95 & & 3.66 & 6.25 & 9.52 & & 0.49 & 0.79 & 1.31 & & 2.88 & 3.56 & 5.03 \\
  \multirow{3}{*}{10}  & 50 & & 13.7 & 23.2 & 29.1 & & 19.2 & 32.5 & 37.5 & & 10.8 & 18.2 & 23.5 & & 13.8 & 23.4 & 32.3 \\
  & 100 & & 8.69 & 15.2 & 23.8 & & 12.4 & 23.6 & 30.9 & & 4.60 & 7.46 & 9.52 & & 9.01 & 13.5 & 18.9 \\
  & 200 & & 2.12 & 4.96 & 7.24 & & 4.02 & 8.17 & 10.6 & & 1.08 & 2.05 & 3.68 & & 3.37 & 6.09 & 10.75 \\
  \multirow{3}{*}{20}  & 50 & & 17.6 & 26.9 & 31.8 & &  28.9 & 37.4 & 39.8 & & 16.5 & 25.3 & 30.0 & & 21.3 & 28.8 & 37.3 \\
  & 100 & & 12.0 & 23.1 & 25.0 & & 16.6 & 30.2 & 30.5 & & 8.21 & 16.2 & 23.9 & & 14.0 & 24.7 & 31.1 \\
  & 200 & & 3.99 & 7.01 & 10.6 & & 5.84 & 10.2 & 15.1 & & 2.49 & 4.05 & 8.46 & & 8.22 & 11.3 & 16.2 \\ \hline
  \multicolumn{5}{l}{RMSE}\\ \hline
  \multirow{3}{*}{5}  & 50 & & 1.78 & 2.63 & 3.88 & & 2.06 & 2.76 & 4.07 & & 0.80 & 0.98 & 1.17 & & 0.88 & 1.04 & 1.21 \\
  & 100 & & 0.87 & 1.04 & 1.19 & & 0.95 & 1.02 & 1.27 & & 0.44 & 0.48 & 0.54 & & 0.50 & 0.52 & 0.57 \\
  & 200 & & 0.69 & 0.64 & 0.70 & & 0.61 & 0.80 & 0.83 & & 0.33 & 0.33 & 0.34 & & 0.42 & 0.39 & 0.41 \\
  \multirow{3}{*}{10}  & 50 & & 4.53 & 7.49 & 9.22 & & 5.50 & 7.79 & 9.29 & & 1.83 & 3.02 & 5.67 & & 2.43 & 3.74 & 6.39 \\
  & 100 & & 1.52 & 1.76 & 2.78 & & 1.59 & 2.00 & 2.41 & & 0.76 & 0.84 & 0.96 & & 0.90 & 0.96 & 1.09 \\
  & 200 & & 0.97 & 1.01 & 1.09 & & 0.94 & 1.21 & 1.24 & & 0.55 & 0.56 & 0.58 & & 0.69 & 0.70 & 0.68 \\
  \multirow{3}{*}{20}  & 50 & & 10.6 & 13.3 & 14.4 & & 11.1 & 13.3 & 14.5 & & 8.48 & 12.8 & 15.1 & & 9.58 & 13.0 & 15.3\\
  & 100 & & 2.61 & 4.06 & 8.25 & & 3.53 & 5.55 & 8.95 & & 1.46 & 1.81 & 2.71 & & 2.21 & 2.67 & 3.98 \\
  & 200 & & 1.46 & 1.65 & 1.78 & & 1.69 & 1.76 & 1.91 & & 0.91 & 0.94 & 1.00 & & 1.13 & 1.22 & 1.31 \\ \hline
\end{tabular}
\end{widepage}
\end{table}

\section{Empirical Analysis}
\subsection{Predicting GDP growth}
We consider the problem of predicting the growth rate of the U.S. quarterly gross domestic product (GDP). In addition, nine macroeconomic variables with different sampling frequencies are
available. The data are obtained from the St. Louis Federal Reserve Economic Data website.
The predictive regression used  is
\begin{eqnarray}
y_{i}=\phi_0+\phi_1 y_{i-1}+\cdots+\phi_a y_{i-a}+\sum_{l=1}^9\sum_{b=0}^{B_l} \beta_{l,b}z_{l,i\times m_l-b}+e_i,
\end{eqnarray}
where $a$ and $B_l$ are nonnegative integers, $y_i$ is the growth rate (first difference of natural  logarithm) of U.S. quarterly seasonally adjusted real GDP, and $z_{l,\cdot}$ are high-frequency covariates with frequency $m_l$, for example, $m_l=3$ for monthly data. The nine covariates considered in this study are as follows: $z_{1,\cdot}$ is the change of monthly civilian unemployment rates; $z_{2,\cdot}$ is the monthly growth rate of all employees' total payrolls; $z_{3,\cdot}$ is the growth rate of the monthly industrial production total index; $z_{4,\cdot}$ is the growth rate of the monthly consumer price index; $z_{5,\cdot}$ is the growth rate of the monthly Moody's Seasoned Baa Corporate Bond Yields; $z_{6,\cdot}$ is the change in the daily 3-Month Treasury Bill Secondary Market Rate; $z_{7,\cdot}$ is the change in the daily 10-Year Treasury Constant Maturity Rate; $z_{8,\cdot}$ is the change in the daily NASDAQ Composite Index; and $z_{9,\cdot}$ is the change in the daily Wilshire 5000 Total Market Full Cap Index. The transformations of all variables are based on those  of \citet{stock2002}. Note that all data are seasonally adjusted, if necessary, and the explanatory variables are monthly or daily data. For daily variables $z_{6,\cdot}$ and $z_{7,\cdot}$, we use data of the first 16 trading days in a month. For daily variables $z_{8,\cdot}$ and $z_{9,\cdot}$, we use data of the first 15 trading days. The sampling period was January 1980 to February 2017, but the prediction origin started with the second quarter of 2013, and ended with the first quarter of 2017.
There was no trading activity during weekends and holidays,
and there exist some missing data in the trading activities.
Trading days for each month vary. We choose the first 15 or 16 trading days,
simply because they are the minimum number of trading days available for each
month (mainly February).

Two types of empirical analysis are examined.
First, we consider a linear model with all explanatory variables,
estimated by the Lasso procedure.
For comparison, we include a model with all explanatory variables except the NASDAQ Composite Index and Wilshire 5000 Total Market Full Cap Index, estimated by the MIDAS regression (denoted by MIDAS-B model), a model with total monthly payrolls for all-employees as the only explanatory variable, also estimated by MIDAS (denoted by MIDAS-A model), and a simple ARMA model of the GDP growth rates (denoted by ARMA model). We use the BIC to select the number of autoregressive lags ($a$) and the lags ($B_l$) of the explanatory variables. The Lasso tuning parameter $\lambda$ is also chosen using the BIC; see \citet{buhlmann2011}. Here, we aggregate the daily explanatory variables $z_6$  and $z_7$ to a weekly frequency for the MIDAS regression.

Table~\ref{tab:table3} shows the median absolute deviation (MAD), mean absolute error (MAE), and root mean squared error (RMSE) for the prediction period. From the table,
it is clear that the Lasso-based model outperforms all the other models in this particular instance.
The poor performance of MIDAS-B is likely due to using too many explanatory variables with multiple
sampling frequencies.

\begin{figure}
\begin{center}
\subfloat[]{\includegraphics[scale=0.47]{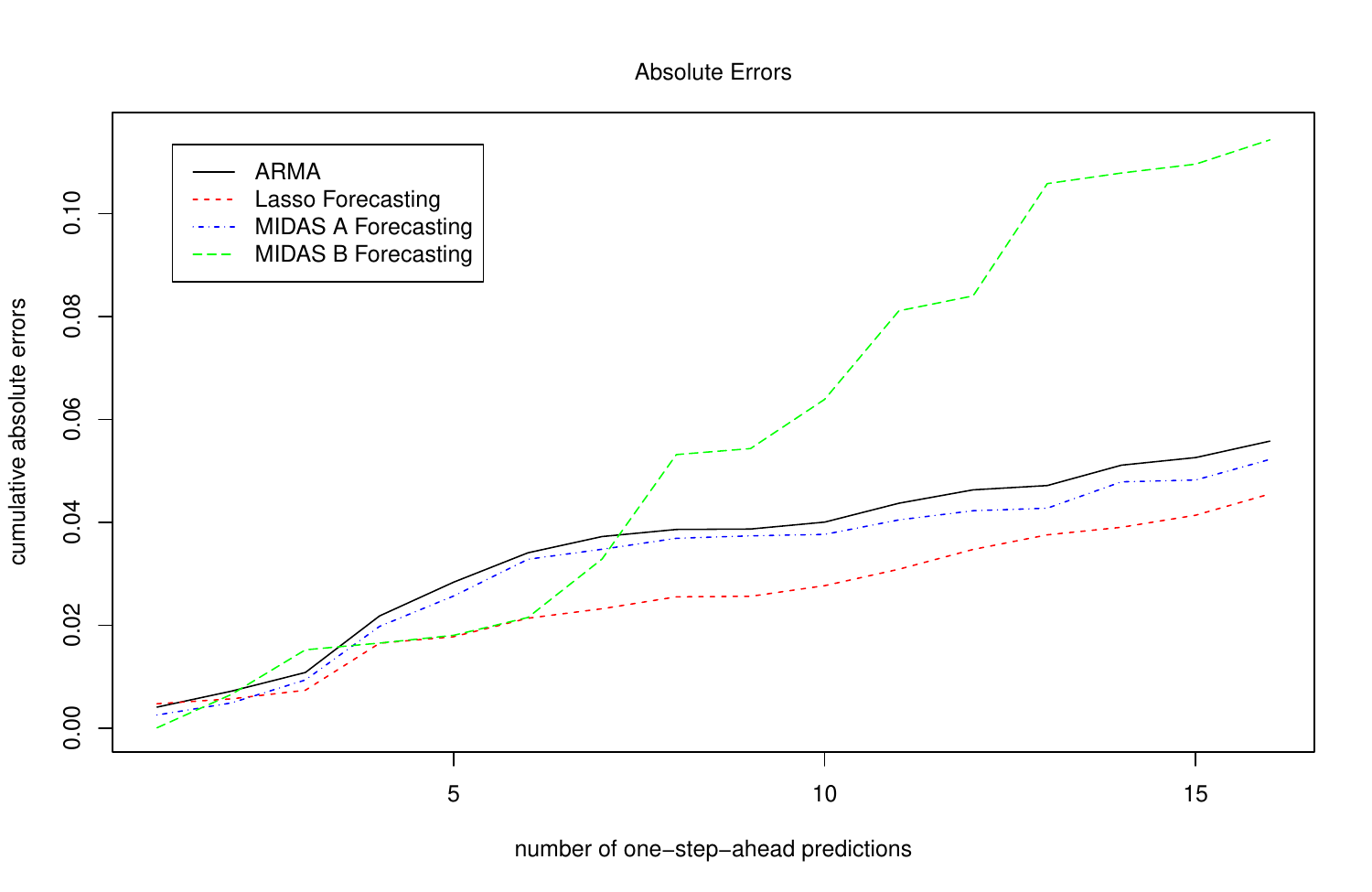}}\\
\subfloat[]{\includegraphics[scale=0.47]{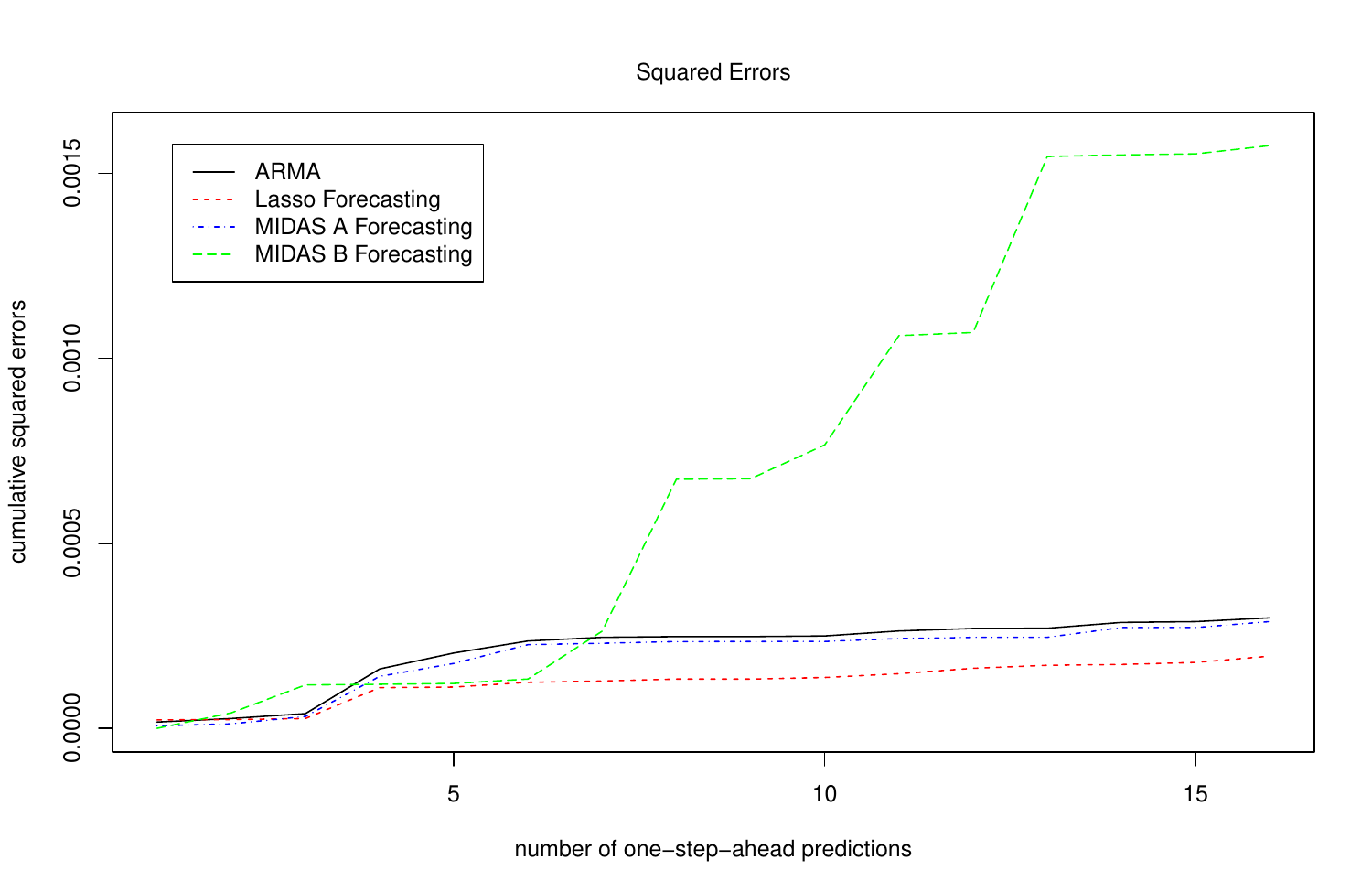}}
\end{center}
\caption{Panel (a): Cumulative absolute errors. Panel (b): Cumulative squared errors. MIDAS-A represents the MIDAS regression model using all-employees' monthly total payrolls as the explanatory variable. MIDAS-B represents the MIDAS regression model with seven regressors $z_{1,\cdot},\cdots,z_{7,\cdot}$, where $z_{6,\cdot}$ and $z_{7,\cdot}$ are aggregated into weekly data.}
\label{fig:fig1}
\end{figure}

Figure \ref{fig:fig1} displays the cumulative absolute errors and the cumulative squared errors for different models in predicting the GDP growth rate. It shows clearly that the Lasso model performs best. The MIDAS-A model also improves the prediction errors over those of the simple ARMA model. However, the MIDAS-B model fares poorly. Consequently, unlike the Lasso model, the MIDAS regression is not robust to the presence of irrelevant regressors. In fact, the MIDAS regression is also sensitive to the weighting schemes and the starting points of its optimization program.

\begin{table}[htbp]
\centering
\caption{Results of out-of-sampling prediction of U.S. quarterly real GDP growth rate.
The data cover the period 1980 to February 2017, but the forecast origins start from
the second quarter of 2013 to the first quarter of 2017. All measurements are  multiplied
by $10^{3}$. In the table, MAD, MAE, and RMSE are the median absolute error, mean
absolute error, and root mean squared error, respectively. }
\label{tab:table3}

\vspace{0.2in}
\begin{tabular}{lcccr}\hline
Model  & MAD & MAE & RMSE \\ \hline
ARMA &3.175 &3.486 &4.319 \\
Lasso  &2.328 &2.845 &3.491 \\
MIDAS-A &2.463 &3.264 &4.245 \\
MIDAS-B &4.089 &7.143 &9.920 \\ \hline
\end{tabular}
\end{table}

Next, we compare between forecasting and nowcasting.
Recall that the goal of nowcasting
is to take advantage of available high-frequency data to
improve the prediction of lower-frequency variables of interest. For the quarterly GDP growth rate,
during the quarter of interest, some monthly macroeconomic variables, and even some daily
economic variables become available. Here, nowcasting attempts to update the GDP prediction
by incorporating the newly available high-frequency explanatory variables. In this exercise,
we consider nowcasting using the first month's data within the quarter, and using the
first two months' data.

For comparison purposes, we employ an autoregressive (AR) model
\begin{equation}\label{ar-forecast}
 y_i = \phi_0 + \phi_1 y_{i-1} + \cdots + \phi_a y_{i-a} +\epsilon_i,
\end{equation}
as a benchmark for prediction. The AR order is selected using the BIC in the modeling subsample,
and is assumed to be fixed in the forecasting subsample. The AR model in Equation (\ref{ar-forecast}) is estimated in two ways. First, it is estimated using the ordinary least squares method,
and we denote the model by AR-OLS. Second, assuming sparsity, we estimate the AR model
via the Lasso method, with the tuning parameter $\lambda$ selected using the BIC. The forecasting result of this model is denoted by AR-Lasso.
These two models represent the performance of forecasting.

For nowcasting, we augment the AR model in Equation (\ref{ar-forecast})
with all explanatory variables available in the first month of the quarter, and denote the
results by Nowcasting 1. Similarly, if we augment the AR model with all
explanatory variables available in the first two months of the quarter, then the
results are denoted by Nowcasting 2.  Specifically, for nowcasting, we employ the model
\[ y_i = \phi_0 + \phi_1 y_{i-1} + \cdots + \phi_a y_{i-a} + \beta^T {\mathbf x}_i + \epsilon_i,\]
where ${\bf x}_i$ denotes the available high-frequency explanatory variables.
For Nowcasting 1, ${\bf x}_i$ consists of data of the first month of a given quarter, whereas
for Nowcasting 2, it consists of data of the first two months of a given quarter. In this
exercise, we use all monthly and daily high-frequency variables $z_{1,\cdot}, \cdots, z_{9,\cdot}$.
We denote the results for the MIDAS regressions as MIDAS-C Nowcasting 1 and MIDAS-C Nowcasting 2,
respectively. Finally, we employ a MIDAS regression that only uses
explanatory variables $z_{1,\cdot}, \cdots, z_{7,\cdot}$ in
the nowcasting and denote the results as MIDAS-D.

Table~\ref{tab:table4} summarizes the performance of nowcasting in predicting U.S. quarterly
GDP growth rates in the forecast period. From the table, we make the following observations.
First, as expected, now-casting fares better than forecasting.
The only exception is MIDAS-D nowcasting.
Second, also as expected, Nowcasting 2 shows some improvement over
Nowcasting 1 for a given model. Keep in mind, however, Nowcasting 1 is available one month
into a quarter, whereas Nowcasting 2 needs to wait for an additional month.
Third, from the performance of MIDAS-C and MIDAS-D,
the stock market indices do not seem to be helpful in predicting the GDP growth rate.
In real applications, there exist many high-frequency explanatory variables, but their contributions
to predicting the low-frequency variable of interest in unknown a priori. In this situation,
our results suggest that the Lasso regression could be helpful.

\begin{figure}
\begin{center}
\subfloat[]{\includegraphics[scale=0.47]{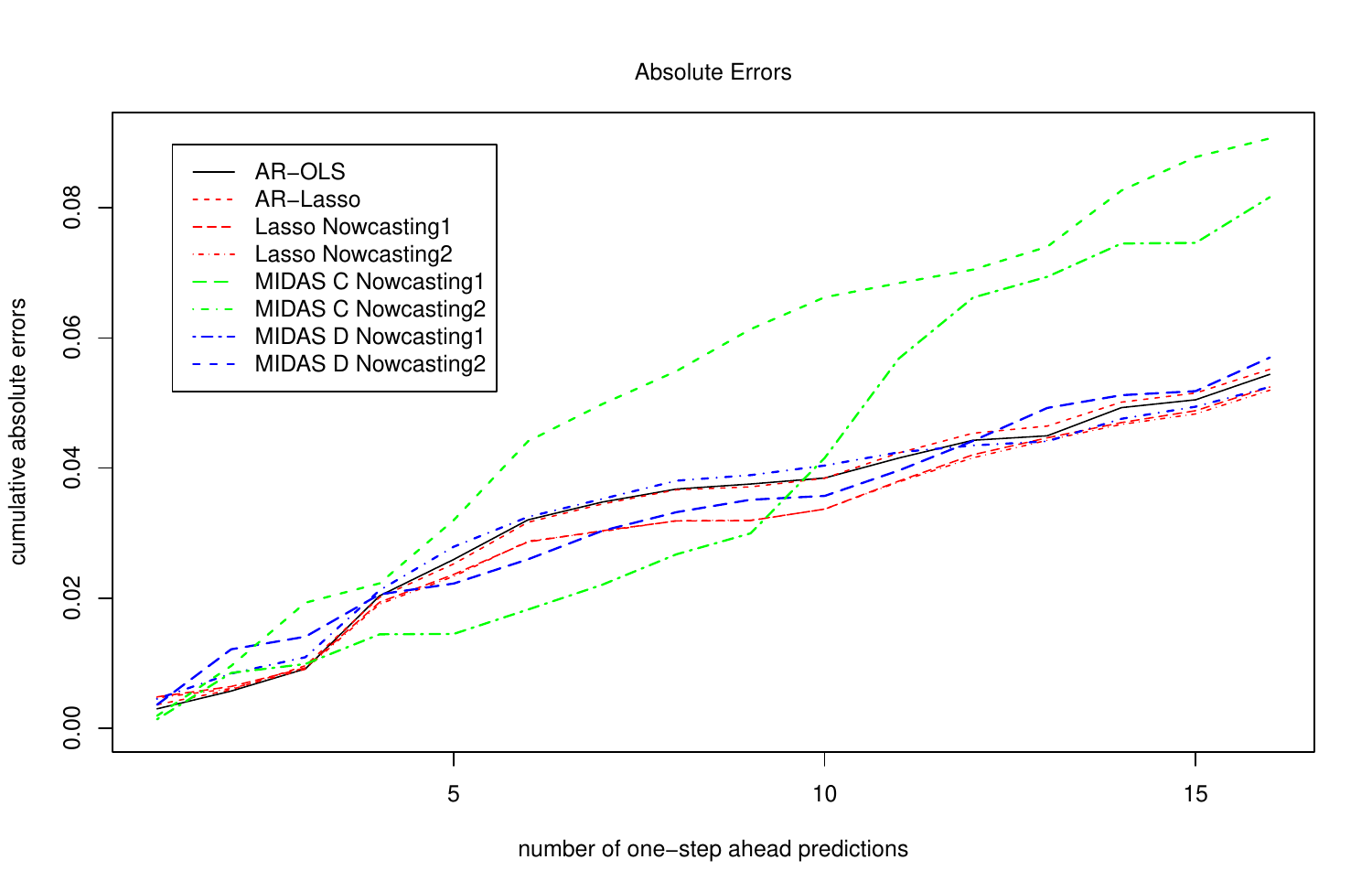}}\\
\subfloat[]{\includegraphics[scale=0.47]{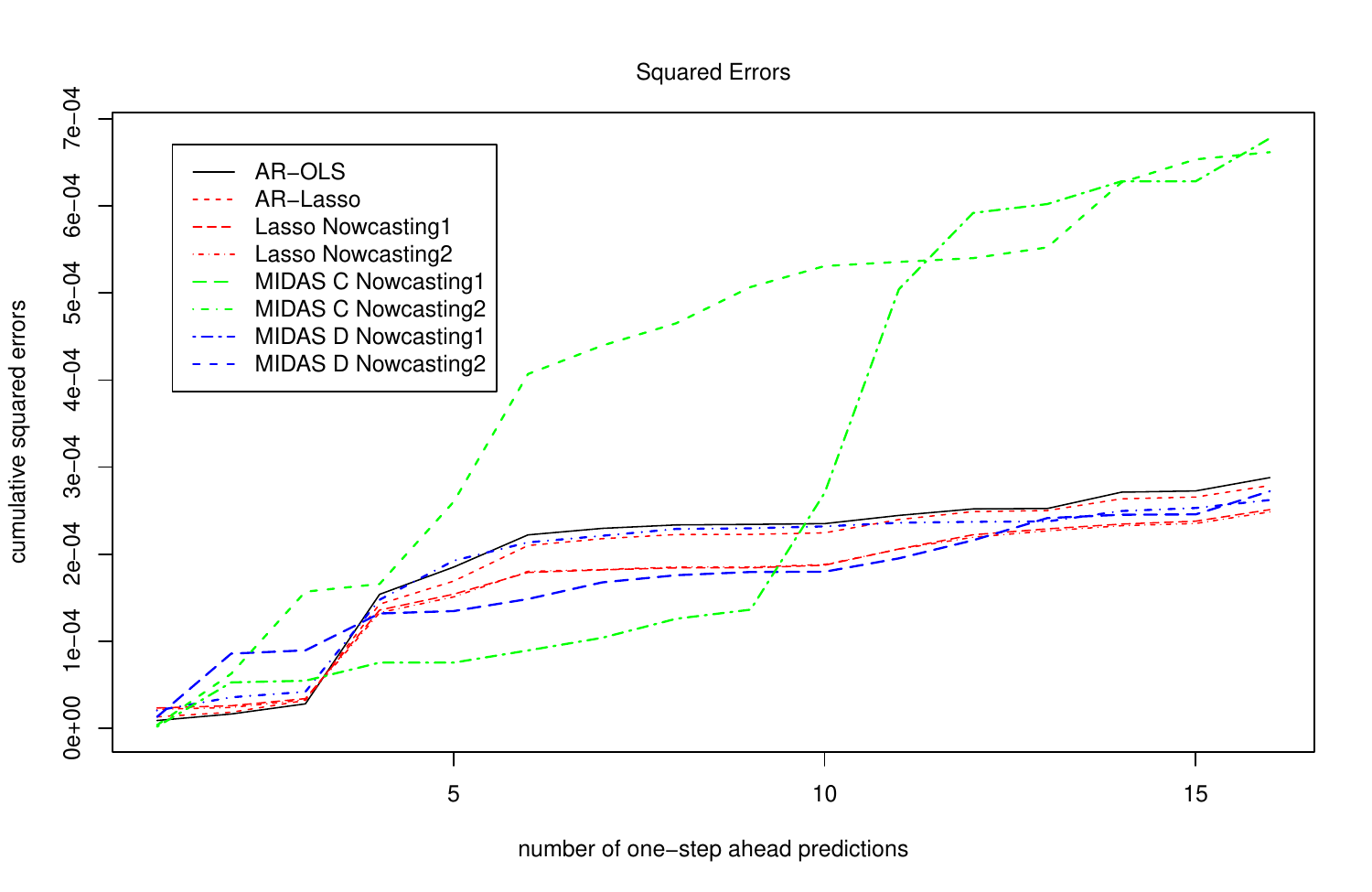}}
\end{center}
\caption{Panel (a): Cumulative absolute errors. Panel (b) Cumulative squared errors. MIDAS-D represents the MIDAS regression model with seven regressors $z_{1,\cdot},...,z_{7,\cdot}$. MIDAS-C represents the MIDAS regression model with nine regressors $z_{1,\cdot},...,z_{9,\cdot}$. Nowcasting 1 and Nowcasting 2 represent predictions of the quarterly GDP growth rate when the first month and the first two months data are available, respectively.}
\label{fig:fig2}
\end{figure}

Figure \ref{fig:fig2} shows that both the Lasso model and the MIDAS-B model improve the prediction via nowcasting. However, when irrelevant variables exist, the MIDAS regression might encounter some
difficulties.

\begin{table}[htbp]
\centering
\caption{Comparison between forecasting and nowcasting in
predicting the U.S. quarterly real GDP growth rate.
The data cover the period 1980 to February 2017, but the forecast origins are from
the second quarter of 2013 to the first quarter of 2017. All measurements are  multiplied
by $10^{3}$. In the table, MAD, MAE, RMSE are the median absolute deviation, mean
absolute error, and root mean squared error, respectively.}
\label{tab:table4}
\begin{tabular}{lccc}\hline
Model  & MAD & MAE & RMSE \\ \hline
AR-OLS &2.865 &3.400 &4.242 \\
AR-Lasso &3.327 &3.448 &4.174 \\
Lasso Now-casting 1 &2.731 &3.278 &3.962 \\
Lasso Now-casting 2 &2.834 &3.247 &3.941 \\
MIDAS-C Now-casting 1 &4.181 &5.102 &6.507 \\
MIDAS-C Now-casting 2 &5.108 &5.666 &6.430 \\
MIDAS-D Now-casting 1 &3.670 &3.561 &4.125 \\
MIDAS-D Now-casting 2 &2.784 &3.279 &4.048 \\ \hline
\end{tabular}
\end{table}

\subsection{Nowcasting PM$_{2.5}$}
Consider next the prediction of PM$_{2.5}$. The response $y$ is the square root transformed daily maximum of PM$_{2.5}$. Hourly data of a monitoring station
in the southern part of Taiwan are used.
To see the nowcasting effects, we consider adding six covariates, which are the
first six hourly PM$_{2.5}$ readings of the same day, starting
from midnight. The sample period is 2006 to 2015, yielding 3650 observations. (Feb 29 was dropped.) We reserve the last 730 data points (two years) for one-step-ahead
out-of-sample forecasts.

For comparison purposes, we first consider the square root PM$_{2.5}$ (i.e., response $y$) as a pure time series. An AR(22) model is selected.
Thus, the baseline model is a univariate AR(22). We denote the model by AR-OLS. For nowcasting, we augment the AR model with the first six hourly readings. If we augment the AR model with the first hourly PM$_{2.5}$ reading, then the results are denoted by Nowcasting 1. Similarly, if we augment the AR model with the first two hourly
PM$_{2.5}$ readings, then the results are denoted by Nowcasting 2, and so on. We denote the results for the autoregressive model with exogenous variables as ARX Nowcasting 1, ARX Nowcasting 2, and so on. We use the BIC to select the number of autoregressive lags. The Lasso tuning parameter $\lambda$ is also chosen using the BIC.

Table \ref{tab:tabled} summarizes the performance of nowcasting in predicting
the daily maximum of PM$_{2.5}$. From the table, we make the following observations. First, as expected, nowcasting outperforms forecasting. Second, also as expected, for a given model, Nowcasting 2 shows some improvement over Nowcasting 1, Nowcasting 3 shows some improvement over Nowcasting 2, and so on. Third, and of most interest, the Lasso estimator significantly outperforms the ARX model and the benchmark model.
In short, the Lasso regression appears helpful in applying nowcasting to PM$_{2.5}$.

\begin{table}[htbp]
\centering
\caption{Comparison between forecasting and nowcasting in
predicting the daily maximum of PM$_{2.5}$.
The data period is 2006 to 2015, and the forecast origins are from
2013 to the end of 2015. (February 29 is excluded). In the table, MAE and RMSE denote the mean
absolute error and root mean squared error for one-step-ahead predictions, respectively.}
\label{tab:tabled}
\begin{tabular}{lccc}\hline
Model  & MAE & RMSE \\ \hline
AR-OLS &1619.6 &73.71 \\
ARX Now-casting 1 &975.9 &46.54 \\
ARX Now-casting 2 &940.9 &44.92 \\
ARX Now-casting 3 &904.2 &43.40 \\
ARX Now-casting 4 &879.7 &42.31 \\
ARX Now-casting 5 &850.6 &41.24 \\
ARX Now-casting 6 &835.3 &40.31 \\
Lasso Now-casting 1 &659.3 &31.74 \\
Lasso Now-casting 2 &628.4 &30.58 \\
Lasso Now-casting 3 &623.2 &30.72 \\
Lasso Now-casting 4 &600.7 &29.63 \\
Lasso Now-casting 5 &595.0 &29.49 \\
Lasso Now-casting 6 &576.3 &28.47 \\\hline
\end{tabular}
\end{table}

\section*{Supplementary Material}
The online Supplementary Material contains proofs of the theorems and lemmas presented in this paper.

\section*{Acknowledgements}
We are grateful to the referees, Associate Editor, and Editor for their many helpful comments.

\nocite{jia2013}
\bibliographystyle{plainnat} 
\bibliography{Lasso}

\end{document}